\documentclass[aap,preprint]{imsart}
\usepackage[utf8]{inputenc}
\usepackage{dsfont}
\RequirePackage[colorlinks]{hyperref}
\usepackage{bbold}
\RequirePackage[numbers]{natbib}

\usepackage{graphicx}
\usepackage{amsmath,amsthm,amssymb}
\usepackage{mathtools,mathdots,mathrsfs}
\usepackage{cleveref}
\usepackage{verbatim}
\usepackage{bm}
\usepackage{subcaption}
\usepackage{natbib}

\usepackage{algorithm, algpseudocode}
\algnewcommand\algorithmicinput{\textbf{Input:}}
\algnewcommand\Input{\item[\algorithmicinput]}
\algnewcommand\algorithmicoutput{\textbf{Output:}}
\algnewcommand\Output{\item[\algorithmicoutput]}

\usepackage[dvipsnames,table]{xcolor}
\usepackage{tikz}
\usetikzlibrary{arrows.meta, bending,     patterns     }

\newtheorem{theorem}{Theorem}[section]
\newtheorem{lemma}[theorem]{Lemma}
\newtheorem{corollary}[theorem]{Corollary}
\newtheorem{proposition}[theorem]{Proposition}
\theoremstyle{remark}
\newtheorem{definition}[theorem]{Definition}
\newtheorem{remark}[theorem]{Remark}
\newtheorem*{example}{Example}

\newcommand\numberthis{\addtocounter{equation}{1}\tag{\theequation}}

\newcommand{\PP}{\mathscr{P}}

\newcommand{\norm}[1]{\ensuremath{\left\| #1\right\|}}

\newcommand{\vect}[1]{\boldsymbol{#1}}
\newcommand{\matr}[1]{\boldsymbol{#1}}

\newcommand{\T}{{\sf T}}        % transposition
\newcommand{\E}{\mathbf{E}} %Expectation
\newcommand{\Proba}{\mathbf{P}}       
\newcommand{\bu}{\vect{u}}

\newcommand{\bM}{\matr{M}}

\newcommand{\bA}{\matr{A}}       
\newcommand{\bQ}{\matr{Q}}       
\newcommand{\bP}{\matr{P}}       

\newcommand{\bB}{\matr{B}}       
\newcommand{\bZ}{\matr{Z}}       
\newcommand{\bU}{\matr{U}}

\newcommand{\bI}{\matr{I}}       
\newcommand{\bV}{\matr{V}}
\newcommand{\bT}{\matr{T}}

\newcommand{\bW}{\matr{W}}       
\newcommand{\bK}{\matr{K}}       
\newcommand{\bR}{\matr{R}}       
\newcommand{\bD}{\matr{D}}       
\newcommand{\bS}{\matr{S}}       
\newcommand{\bL}{\matr{L}}
\newcommand{\R}{\mathbb{R}}       
\newcommand{\tbL}{\widetilde{\matr{L}}}
\newcommand{\tbU}{\widetilde{\matr{U}}}
\newcommand{\tbLam}{\widetilde{\bm{\Lambda}}}

\newcommand{\X}{\mathcal{X}}
\newcommand{\Xe}{\mathcal{X}_\varepsilon}
\newcommand{\Xs}{\mathcal{X}_\star}

\newcommand{\Y}{\mathcal{Y}}

\newcommand{\flatlim}{\varepsilon\rightarrow0}

\renewcommand{\O}{\mathcal{O}}       
\DeclareMathOperator{\rank}{rank}
\DeclareMathOperator{\mspan}{span}
\DeclareMathOperator{\orth}{orth}

\DeclareMathOperator{\diag}{diag}
\DeclareMathOperator{\Tr}{Tr}

\definecolor{darkgreen}{rgb}{0,0.6,0}
\newcommand{\fin}{\color{black}}

\newcommand{\ELE}[2]{\begin{pmatrix} #1 ; #2  \end{pmatrix}}

\newcommand{\ones}{\vect{\mathbf{1}}}

\newcommand{\Ind}{\vect{\mathds{I}}}

\newcommand{\mDPP}[1]{|DPP|_{#1}}

\newcommand{\ppDPP}{DPP}
\newcommand{\mppDPP}[1]{|DPP|_{#1}}

\newcommand{\cW}{\mathcal{W}}
\newcommand{\cA}{\mathcal{A}}

\newcommand{\Qort}{\matr{Q}_{\bot}}

\newcommand{\gL}{\vect{{\mathcal{L}}}}
\begin{document}
\begin{frontmatter}

\title{Extended L-ensembles: a new representation for Determinantal Point Processes}
\runtitle{Extended L-ensembles}

\begin{aug}
  \author[A]{\fnms{Nicolas}
    \snm{Tremblay}\ead[label=e2]{nicolas.tremblay@gipsa-lab.fr}},
  \author[A]{\fnms{Simon} \snm{Barthelm\'e}\ead[label=e1]{name.surname@gipsa-lab.fr}},
  \author[B]{\fnms{Konstantin} \snm{Usevich}\ead[label=e3]{konstantin.usevich@univ-lorraine.fr}}
  \and
  \author[A]{\fnms{Pierre-Olivier} \snm{Amblard}\ead[label=e4]{pierre-olivier.amblard@gipsa-lab.fr}}

\address[A]{CNRS, Univ. Grenoble Alpes,  Grenoble INP, GIPSA-lab}%\ \printead{e1}}

  \address[B]{Universit\'{e} de Lorraine and CNRS, CRAN (Centre de Recherche en
    Automatique de Nancy)}% \printead{e3}}
\end{aug}

  \begin{abstract}
    Determinantal point processes (DPPs) are a class of repulsive point
    processes, popular for their relative simplicity. They are traditionally
    defined via their marginal distributions, but a subset of DPPs called
    ``L-ensembles'' have tractable likelihoods and are thus particularly easy to
    work with. Indeed, in many applications, DPPs are more naturally defined
    based on the L-ensemble formulation rather than through the marginal kernel.
    
    The fact that not all DPPs are L-ensembles is unfortunate, but there is a
    unifying description. We introduce here \emph{extended} L-ensembles, and
    show that all DPPs are extended L-ensembles (and vice-versa). Extended
    L-ensembles have very simple likelihood functions, contain L-ensembles and
    projection DPPs as special cases. From a theoretical standpoint, they fix
    some pathologies in the usual formalism of DPPs, for instance the fact that
    projection DPPs are not L-ensembles. From a practical standpoint, they
    extend the set of kernel functions that may be used to define DPPs: we show
    that conditional positive definite kernels are good candidates for defining
    DPPs, including DPPs that need no spatial scale parameter.

    Finally, extended
    L-ensembles are based on so-called ``saddle-point matrices'', and we prove
    an extension of the Cauchy-Binet theorem for such matrices that may be of
    independent interest.
  \end{abstract}
\end{frontmatter}

\pagebreak
\tableofcontents
\pagebreak

\section*{Introduction}
\label{sec:intro}
Determinantal point processes are by now perhaps the most famous example of
repulsive point processes. They appear as a model for the position of
fermionic particles in an energy potential \cite{Macchi:CoincidenceApproach},
but also occur in random matrix theory and graph theory, for instance as the
distribution of edges in uniform spanning trees \cite{burton1993local}. DPPs can
also be used as tools designed for particular computational tasks: for instance,
they have been advocated in machine learning as a way of providing samples with
guaranteed diversity \cite{kulesza2012determinantal}, as a way of performing
Monte Carlo integration \cite{bardenet2016monte,coeurjolly2020monte}, or 
accelerating classical linear algebra tasks such as regression or low rank approximation of matrices~\cite{belhadji_determinantal_2020,derezinski_determinantal_2021, tremblay2019determinantal}.

We will use discrete sampling problems as a motivating example throughout.  In that framework, one has a set of
$n$ items, and one desires to produce a subset $\X$ of size $m \ll n$ such that
no two items in $\X$ are excessively similar.  A key aspect of DPPs is that
``diversity'' is defined relative to a notion of similarity represented by a
positive-definite kernel. For instance, if the items are vectors in $\R^d$,
similarity may be defined via the squared-exponential (Gaussian) kernel:
\begin{equation}
  \label{eq:squared-exp-kernel}
  \kappa(\vect{x},\vect{y}) = \exp \left( -  \norm{\vect{x}-\vect{y}}^2 \right)
\end{equation}
Here $\vect{x}$ and $\vect{y}$ are two items, and similarity is a decreasing
function of distance.

%  a large subclass called \emph{L-ensembles} grouping the DPPs that can sample the empty set (the probability of sampling the empty set is strictly positive); and a much smaller class grouping DPPs that cannot (the probability is strictly zero). Precise definitions are to be found in section~\ref{sec:definitions}. 

The class of DPPs can be separated into two subclasses: L-ensembles and the
rest. By definition, an L-ensemble based on the $n\times n$ kernel matrix $\bL=[\kappa(\vect{x}_i,\vect{x}_j)]_{i,j}$ is a distribution over random subsets $\X$ such that:
\[ \Proba(\X = X) \propto \det \bL_{X} \]
where $\bL_{X}$ is the principal submatrix of $\bL$ indexed by $X$. 
If two or more points in $X$ are very similar (in the sense of the kernel function), then the matrix $\bL_X$
has rows that are nearly collinear and the determinant is small. This in turns
makes it unlikely that such a set $\X$ will be selected by the L-ensemble. 

L-ensembles are thus highly intuitive, and it is easy to design a DPP
appropriate for a particular situation, since one only needs to pick an
appropriate kernel. Unfortunately, not all DPPs are L-ensembles. This includes
some naturally-occuring DPPs, like the uniform spanning tree, or the roots of
uniform spanning forests \cite{avena_two_2017}. In this work we introduce
extended L-ensembles, which are easy to work with, and include L-ensembles as a
special case. We show that all DPPs are extended L-ensembles, and vice versa.

Extended L-ensembles are based on an $n \times n$ kernel matrix, noted $\bL$,
augmented by a set of $p$ vectors, noted $\bV \in \R^{n \times p}$. $\bV$
represents a set of ``obligatory'' features, as we explain below, and the
likelihood function is just:
\[ \Proba(\X=X) \propto  \det \begin{pmatrix}
    {\bL}_{X} & \bV_{X,:} \\
    (\bV_{X,:})^\top & \matr{0} 
	\end{pmatrix}\]
	where $\bV_{X,:}$ is the matrix $\bV$ restricted to its lines indexed by $X$.

As a practical benefit, extended L-ensembles broaden the class of kernels that
may be used to define DPPs, to include \emph{conditional positive-definite
  kernels}. This allows for a class of ``default'' DPPs that do not require a
length-scale parameter\footnote{Since this work first appeared in our prior work in arXiv:2007.04117,  extended
  L-ensembles have been used in \cite{fanuel2020determinantal} to define Nyström
approximations appropriate for semi-parametric models.}.

\subsection*{Structure of the paper}
We begin with some definitions and background in section \ref{sec:definitions}.
Section \ref{sec:ppDPPs} defines extended L-ensembles while
section~\ref{sec:general_properties} gives some of their major properties.
As a theoretical application, section~\ref{sec:ppDPP-as-limits} shows how
extended L-ensembles arise in perturbative limits of DPPs. As a practical
application, section~\ref{sec:ppdpp-examples-cpdef} shows that some interesting extended L-ensembles can be
constructed via
conditionally (semi)-definite positive functions.

We chose to focus on discrete DPPs, for simplicity. The continuous case is a
straightforward generalisation of our results, which we sketch in 
our concluding remarks (section~\ref{sec:conclusion}).

\section{Definitions and background}
\label{sec:definitions}

We briefly recall some  definitions on DPPs along with
fixed-size DPPs, a useful variant (as well as L-ensembles and fixed-size L-ensembles). For details we refer the reader to
\cite{Barthelme:AsEqFixedSizeDPP} and  \cite{KuleszaTaskar:FixedSizeDPPs}. All
of the results below are classical.

\subsection{Some determinant lemmas}
\label{sec:determinant-lemmas}
We start with notations and a few useful determinant lemmas.

 Let $\matr{A}$ be a $n \times n$ matrix, and
$Y$, $Z$ be two subsets of indices. Then $\matr{A}_{Y,Z}$ is the submatrix
of $\matr{A}$ formed by retaining the rows in $Y$ and the columns in $Z$.  Furthermore, $\matr{A}_{:,Y}$ (resp. $\matr{A}_{Y,:}$) is the matrix made of the full columns (resp. rows) indexed by $Y$. Finally, we
let $\matr{A}_Y = \matr{A}_{Y,Y}$.
Also, for a matrix $\matr{V}$, by $\mspan(\matr{V})$ we denote its column span, and by $\orth(\matr{V})$ the orthogonal complement of $\mspan(\bV)$. 

We shall need a number of basic results on determinants. The Cauchy-Binet lemma is central to the theory of DPPs and generalises the
well-known relationship $\det(\matr{A}\matr{B})=\det(\matr{A})\det(\matr{B})$
(for square $\matr{A}$ and $\matr{B}$) to rectangular matrices. 

\begin{lemma}[Cauchy-Binet]
  \label{lem:cauchy-binet}
  Let $\matr{M} = \matr{A} \matr{B}$, with $\matr{A}$ a $m \times n$ matrix and  $\matr{B}$ a $n \times m$ matrix ($n\geq m$). Then:
  \begin{equation}
    \label{eq:CauchyBinet}
    \det \matr{M} = \sum_{Y, |Y| = m} \det \matr{A}_{:,Y} \det \matr{B}_{Y,:}
  \end{equation}
  where the sum is over all subsets $Y \subseteq \{1, \ldots, n\}$ of size $m$.
\end{lemma}
We will also frequently use the following simple corollary of the Cauchy-Binet lemma.
\begin{corollary}
Let $\matr{M} =  \matr{U} \matr{\Lambda} \matr{U}^{\top}$, where $\matr{U}$ is $m\times n$ , $n \ge m$ and $\matr{\Lambda}$ is a diagonal matrix.
Then:
\[
\det \matr{M} = \sum_{Y, |Y| = m} (\det (\matr{U}_{:,Y}))^2 \det(\matr{\Lambda}_{Y}).
\]
\end{corollary}
The next result is a well-known determinantal counterpart of the
Sherman-Woodbury-Morrisson lemma:

\begin{lemma}
  \label{lem:low-rank-update-det}
  Let $\matr{A}$ be an invertible matrix of size $n \times n$, $\matr{U}$
  of size $n \times m$, and $\matr{W}$ an invertible matrix of size $m \times m$. Then it holds that:
  \begin{equation}
    \label{eq:determinant-update}
    \det(\matr{A}+\matr{UWU}^\top)=\det(\matr{A})\det(\matr{W})\det(\bW^{-1}+\matr{U}^\top\matr{A}^{-1}\matr{U}).
  \end{equation}
\end{lemma}

Finally, a related lemma is useful for block matrices:

\begin{lemma}
  \label{lem:block-det}
  Let $\matr{M} =
  \begin{pmatrix}
    \matr{A} & \matr{U} \\
    \matr{U}^\top & \matr{W}
  \end{pmatrix}
  $,  with $\matr{A}$ invertible.  Then
  \begin{equation}
    \label{eq:block-det}
    \det(\matr{M}) = \det(\matr{A})\det(\matr{W}-\matr{U}^\top\matr{A}^{-1}\matr{U}).
  \end{equation}
\end{lemma}

The next two lemmas concern so-called ``saddle-point matrices'', and are proved in
\cite[Appendix A]{BarthelmeUsevich:KernelsFlatLimit}.  

\begin{lemma}[{\cite[Lemma 3.10]{BarthelmeUsevich:KernelsFlatLimit}}]
  \label{lem:det-saddlepoint}
  Let $\bL \in \R^{n \times n}, \bV \in \R^{n \times p}$, 
with $\bV$ of full column rank and $p \leq n$. 
  Let 
  $\bQ \in \R^{n \times (n-p)}$ be an  orthonormal basis for $\orth({\bV})$ (i.e., $\bQ^{\top} \matr{V} = \matr{0}$, $\rank(\bQ) = n-p$).  
  Then:
  \begin{equation}
    \label{eq:det-saddlepoint}
  \det       \begin{pmatrix}
    \bL & \bV \\
    \bV^\top & \matr{0} 
  \end{pmatrix} = (-1)^p \det(\bV^\top\bV)\det(\bQ^\top\bL\bQ).
\end{equation}
\end{lemma}

In the next lemma, we use $[t^r] g(t)$ to denote the coefficient corresponding to $t^r$ in the power series $g$. 
For instance, if $g(t) = 1-t^2+2t^3$, then $[t^0]g(t) = 1$ and $[t^3]g(t) = 2$. 
\begin{lemma}[{\cite[Lemma 3.11]{BarthelmeUsevich:KernelsFlatLimit}}]\label{lem:det-coef-polynomial}
Let $\bL \in \R^{n \times n}$ and $\bV \in \R^{n\times p}$.
Then:
\[
[t^p]\det( \bL +t \bV \bV^\top ) =  (-1)^{p} \det\begin{pmatrix}
      \matr{L} &  \matr{V} \\
      \matr{V}^{\T} & 0 
    \end{pmatrix}. 
\] 
\end{lemma}

\begin{remark}\label{rem:det-coef-degree}
The polynomial  $g(t) = \det( \bL +t \bV \bV^\top )$ is of degree at most $p$, i.e., lemma~\ref{lem:det-coef-polynomial} gives the coefficient for the highest possible power of $t$.
While this remark is missing in the original statement of
lemma~\ref{lem:det-coef-polynomial} (see \cite[Lemma
3.11]{BarthelmeUsevich:KernelsFlatLimit}), it can be easily verified by
inspecting the proof of the lemma in \cite[Appendix
A]{BarthelmeUsevich:KernelsFlatLimit}.
\end{remark}

\subsection{Determinantal processes}
\label{sec:def-DPPs}
\subsubsection{DPPs}
Let $\Omega=\{ \vect{x}_1, \ldots, \vect{x}_n\} \subset \R^d$ be a collection of
vectors called the \emph{ground set}. A finite point process $\X$ is a random
subset $\X \subseteq \Omega$. Abusing notation, we sometimes use $\X$ to
designate the indices of the items, rather than the items themselves. Which one we mean should be clear from context.

\begin{definition}[Determinantal Point Process]\label{def:dpp}
	Let $\bK \in \R^{n \times n}$ be
	a positive semi-definite matrix verifying $\bm{0}\preceq \bK \preceq \bm{I}$. $\X$ is a DPP with
	marginal kernel $\bK$ if
	\begin{equation}
	\label{eq:def-dpp}
	\forall A\subseteq\Omega\qquad\Proba(A\subseteq\mathcal{X}) = \det \bK_A,
	\end{equation}
	where by convention, $\det \bK_{\varnothing} = 1$. 
\end{definition}
This definition is the historical one \cite{Macchi:CoincidenceApproach} and
determines what we will refer to as the \emph{class of DPPs}. However,
manipulating inclusion probabilities rather than the joint probability
distribution itself is often cumbersome. This usually leads authors to consider
a slightly less general class of DPPs: the L-ensembles \cite{borodin2005eynard}. 
\begin{definition}[L-ensemble]\label{def:dpp_via_L}
	Let $\bL \in \R^{n \times n}$ designate
	a positive semi-definite matrix. An L-ensemble based on $\bL$ is a point process $\X$ defined as
	\begin{equation}
	\label{eq:def-dpp_via_L}
	\Proba(\X=X) = \frac{\det \bL_X}{Z},
	\end{equation}
	where by convention, $\det \bL_{\varnothing} = 1$. Thus: $\Proba(\X=\varnothing) =1/Z>0$. 
\end{definition}
In Eq. \eqref{eq:def-dpp_via_L}, \fin $Z = \sum_{X \subseteq \Omega} \det \bL_X$ is a normalisation constant and can be shown~\cite{kulesza2012determinantal} to equal
$\det(\bI + \bL)$. 

L-ensembles are indeed a subclass of DPPs:
\begin{lemma}\label{lem:LisDPP}
	An L-ensemble based on the positive semi-definite matrix $\bL$ is a DPP. It is noted  $\X\sim DPP(\bL)$ and its marginal kernel verifies
	\begin{align}
	\label{eq:connection_KL}
	\bK=\bL(\bI+\bL)^{-1}.
	\end{align}
\end{lemma}
\begin{proof}
	See, \emph{e.g.}, Thm 2.2 of~\cite{kulesza2012determinantal}; or the discussion in  Appendix~\ref{sec:marginal-kernel-ppDPPs-proof}.
\end{proof}

L-ensembles are in fact a strict subset of all DPPs:
\begin{lemma}
	\label{lem:DPPsubclass}
	A DPP with marginal kernel $\bK$ is an L-ensemble if and only if $\bK$ verifies $\bm{0}\preceq\bK\prec\bI$ (note the $\prec$ sign, implying that no eigenvalue of $\bK$ is allowed to be equal to one). If $\X$ is a DPP with such a marginal kernel, then $\X\sim DPP(\bL)$, with $\bL$ verifying:
	$$\bL=\bK(\bI-\bK)^{-1}.$$
\end{lemma}
\begin{proof}
	($\Leftarrow$) If $\bK$ does not contain any eigenvalue equal to 1, then Eq.~\eqref{eq:connection_KL} inverts as $\bL=\bK(\bI-\bK)^{-1}.$ ($\Rightarrow$) We show the contraposition. If $\X$ is a DPP with a marginal kernel $\bK$ containing at least one eigenvalue equal to one, then its size $|\X|$ is necessarily larger than one (see lemma~\ref{lemma:bernoulli}). Thus, it cannot be an L-ensemble (L-ensembles have a non-null probability of sampling $\varnothing$). 
\end{proof}
\begin{remark}
As a consequence, the class of DPPs can be separated in two: the L-ensembles
(all DPPs with marginal kernel verifying $\bm{0}\preceq \bK\prec\bI$), and the rest (all DPPs with marginal kernel whose spectrum contains at least one
eigenvalue equal to one). \end{remark}

In DPPs, the size (cardinal) of $\X$, denoted by $|\X|$, is a random variable. Its
distribution is as follows \cite{Hough:DPPandIndep}:

\begin{lemma}
	\label{lemma:bernoulli}
	Let $\bm{0}\preceq\bK\preceq\bI$ be a marginal kernel with eigenvalues $\mu_1,\ldots, \mu_{n}$. Let $\X$ be a DPP with this marginal kernel. Then, $|\X|$ has the same distribution as $\sum_{i=1}^n B_i$, 
	where $B_i$ is a Bernoulli random variable with expectation
	$\E(B_i)=\mu_i$, and the $B_i$'s are distributed independently.
	In particular, the expected size of the DPP, $\E(|\X|)$, can be directly
	deduced from the above to be
  \begin{equation}
    \label{eq:trace-K}
    \E(|\X|)= \sum \mu_i = \Tr (\bK)
  \end{equation}
\end{lemma}
\subsubsection{Fixed-size DPPs}
The cardinal of a DPP is thus in general random. Such varying-sized samples are
not practical in many applications (one desires a subset of size 50, not
a subset of size 50 on average but which may be of size 35 or 56); which led the authors of~\cite{kulesza2011k} to define fixed-size DPPs\footnote{They are often called k-DPPs in the literature, but we prefer ``fixed-size DPPs'' in order not to overload the
	symbol $k$ too much.}
\begin{definition}[Fixed-size Determinantal Point Process]
	\label{def:fsDPP}
	A fixed size DPP of size $m$ is a DPP $\X$ conditioned on $|\X|=m$. 
\end{definition}
A subclass of fixed-size DPPs is the class of fixed-size L-ensembles: 
\begin{definition}[Fixed-size L-ensemble]
	\label{lemma:mDPP_via_L-ens}
  Let $\bm{0}\preceq\bL$ be a positive semi-definite matrix. A fixed-size L-ensemble is a point process $\X$ defined as:
  \begin{equation}
    \label{eq:def-fsdpp}
  \Proba(\X=X) =    \begin{cases}
 \displaystyle     \frac{\det \bL_{X}}{Z_m} & \mathrm{if} ~|X|=m,  \\
      0 & \mathrm{otherwise.}
    \end{cases} 
  \end{equation}
  where $Z_m$ is the normalisation constant. 
\end{definition}
Using the indicator function $\Ind(\cdot)$, we may rewrite Eq.~\eqref{eq:def-dpp} more compactly as:
\[  \Proba(\X=X) = \frac{\det \bL_{X}}{Z_m} \Ind(|X|=m).\]
\begin{lemma}
	A fixed-size L-ensemble is a fixed-size DPP, and we write it $\X \sim \mDPP{m}(\bL)$. 
\end{lemma}

We use the notation $\X \sim \mDPP{m}(\bL)$ to distinguish from (standard)
random-size L-ensembles.

It is important to understand that, in general, fixed-size DPPs are not DPPs, with the exception of projection DPPs (see Sec.~\ref{sec:diag-and-proj}). In particular, whereas all DPPs have a marginal kernel, fixed-size DPPs (again with the exception of projection DPPs) do not have marginal kernels: there does not exist a matrix whose principal minors are the marginal probabilities. The question of inclusion probabilities in fixed-size
DPPs is treated at length in \cite{Barthelme:AsEqFixedSizeDPP}.

The constant $Z_m = \sum_{\X,|\X|=m} \det \bL_\X$ in Eq.~\ref{eq:def-dpp} is a normalisation constant and one can show that it 
equals the $m$-th ``elementary symmetric polynomial'' of $\bL$, a quantity that
depends only on the spectrum of $\bL$, and plays an important role in the
theory of DPPs. 
\begin{lemma}[{\cite[Theorem 1.2.12]{horn1990matrix}}]
  \label{lem:esp}
  Let $\bL \in \R^{n \times n}$ be a  matrix with eigenvalues $\lambda_1,\ldots,\lambda_n$. The
  $m$-th elementary symmetric polynomial of $\bL$ is defined as:
  \begin{equation}
    \label{eq:esp}
    e_m(\bL) := \sum_{|X|=m} \prod_{i \in X} \lambda_i,
  \end{equation}
  i.e., $e_0(\bL) = 1$, $e_1(\bL) = \sum_{i } \lambda_i = \Tr(\bL) $, $e_2(\bL) = \sum_{i    < j} \lambda_i  \lambda_j, \; \ldots, e_n(\bL) = \det(\bL)$. 
  One has:
   \begin{equation} 
    \label{eq:esp-det}
  Z_m =   \sum_{X,|X|=m} \det \bL_X = e_m(\bL).
  \end{equation}
\end{lemma}
Since $e_m(\bL)$ is the sum of all the principal minors of fixed size $m$, we immediately
obtain the following corollary on the distribution of the size of an L-ensemble:
\begin{corollary}
  The probability that $\X \sim DPP(\bL)$ has size $m$ is given by:
  \begin{equation}
    \label{eq:distr-size-dpp}
    p(|\X|=m) = \frac{e_m(\bL)}{e_0(\bL)+e_1(\bL)+\ldots+e_n(\bL)}.
  \end{equation}
\end{corollary}

\begin{remark}\label{rm:mixture}
  Since a fixed-size L-ensemble is just an L-ensemble conditioned on size, an L-ensemble may also be viewed as a mixture of fixed-size L-ensembles. The size $m$ can be drawn according to its marginal
  distribution (Eq.~\eqref{eq:distr-size-dpp}), and conditional on $|\X|=m$, the fixed-size L-ensemble can be sampled.
\end{remark}

\subsubsection{Two useful special cases}
\label{sec:diag-and-proj}

There are two special cases of (fixed-size) DPPs that are useful to
study on their own, both from a practical and theoretical viewpoint. These are
the DPPs with \emph{diagonal} kernels and those with \emph{projection} kernels. 

As it will be shown in section~\ref{sec:mixture-representation}, these two examples are the key components for sampling any DPP using the mixture representation.

\paragraph{Diagonal kernels} Diagonal L-ensembles are in a way the most basic kind of DPPs 
(although the fixed-size case is surprisingly intricate).

\begin{lemma}
  \label{lem:diagonal-dpp-bernoulli}
  An L-ensemble based on a diagonal positive semi-definite matrix $\bL$, $\Y \sim DPP(\matr{\Lambda})$ with $\matr{\Lambda}=\diag(\lambda_1\ldots,\lambda_n)$, 
  is a Bernoulli process: each event $i \in \Y$ is independent and occurs with
  probability $\pi_i=\frac{\lambda_i}{1+\lambda_i}$.
\end{lemma}
\begin{proof}
  \begin{align*}
    \Proba(\Y = Y) &= \frac{\prod_{i\in Y} \lambda_i}{\det(\bI + \matr{\Lambda})} 
               = \frac{\prod_{i\in Y} \lambda_i}{\prod_{j=1}^n (1+\lambda_j)} = \left(\prod_{i\in Y} \frac{\lambda_i}{1+\lambda_i}\right)\left(\prod_{j \in Y^c}  \frac{1}{1+\lambda_i}\right)\\
               &= \prod_{i=1}^n (\pi_i)^{B_i}(1-\pi_i)^{B_i}
  \end{align*}
  where $B_i$ is the Bernoulli variable indicating $i \in \Y$. 
\end{proof}
\begin{remark}
  For fixed-size L-ensembles this is no longer true: $\Y \sim
  \mDPP{m}(\matr{\Lambda})$ is not a Bernoulli process, as the
  events are no longer independent but indeed negatively associated. To see why,
  note that since the total size is fixed, conditional on $i \in \Y$ other
  points are less likely to be included. 
\end{remark}
\begin{remark}
  $\Y \sim \mDPP{m}(\bI)$ is a uniform sample of size $m$ without replacement. 
\end{remark}
Fixed-size diagonal L-ensembles have been studied at some length in the past, notably in the sampling survey literature. Many important features of these processes were reported in \cite{ChenDL94}.

\paragraph{Projection DPPs}

Projection DPPs designate DPPs formed from projection matrices.
Projection DPPs have many unique features, for instance that of being both DPPs
and fixed-size DPPs. Section \ref{sec:ppDPPs} will introduce a
generalisation called ``partial projection DPPs''. 
The definition of a projection DPP is as follows: 

\begin{definition}[Projection DPP]
	Let $\bU$ be an $n \times m$ matrix with $\bU^\top \bU = \bI_m$. A projection DPP is a DPP with marginal kernel $\bK=\bU\bU^\top$.
\end{definition}

The name ``projection DPP'' comes from the fact that $\bU \bU^\top$ is a projection
matrix (its eigenvalues are 1, with multiplicity $m$, and 0 with multiplicity
$n-m$). As $\bK$'s spectrum contains at least an eigenvalue equal to 1, a projection DPP is \emph{not} an L-ensemble (see lemma~\ref{lem:DPPsubclass}). However, a projection DPP can be equivalently defined as a fixed-size L-ensemble:

\begin{lemma}[See e.g., {\cite[Lemma 1.3]{Barthelme:AsEqFixedSizeDPP}}]
  \label{lem:marginal-kernel-proj}
  Let $\bU$ be an $n \times m$ matrix with $\bU^\top \bU = \bI_m$. A projection DPP with marginal kernel $\bU \bU^\top$ is a fixed-size L-ensemble $\X \sim \mDPP{m}(\bU\bU^\top)$.
\end{lemma}

In fact, the only class of fixed-size DPPs that admit a marginal kernel are the projection DPPs. 
The next result states that a projection DPP is what one obtains when
sampling a fixed-size L-ensemble of size $m$ from a positive semi-definite matrix $\bL$ of rank $m \leq n$. 
\begin{lemma}
  \label{lem:max-rank-dpp}
  Let $\X \sim \mDPP{m}(\bL)$, with $\rank(\bL) = m$, and let $\bU\in \R^{n \times m}$ denote an
  orthonormal basis for $\mspan \bL$.  Then, equivalently, $\X \sim
  \mDPP{m}(\bU \bU^\top)$
\end{lemma}
\begin{proof}
  Given the assumptions, we may write $\bL = \bU \bM \bM^\top \bU^\top$ with $\bU \in
  \R^{n \times m}$, and $\bM \in \R^{m \times m}$.
  Now, bearing in mind that $|\X| = m$, we have:
  \begin{align*}
    \Proba(\X = X) &\propto \det(\bL_X)  = \det(\bU_{X,:} \bM )^2 \propto  \det(\bU\bU^\top)_X
  \end{align*}
  where we used the fact that $\bU_{X,:}$ is square and that $\det(\bM)$
  is independent of $X$. 
  Note that any orthonormal basis works, for instance the eigenvectors of $\bL$ associated with a non-null eigenvalue,
  but not only: the Q factor in the QR factorisation of $\bL$ would work as well.
\end{proof}

\begin{remark}
Note that lemma~\ref{lem:max-rank-dpp} is valid only for fixed-size L-ensembles with rank of $\bL$ exactly equal to $m$.
In the case $\rank \bL > m$, the fixed-size L-ensemble $\X \sim \mDPP{m}(\bL)$ is no longer a projection DPP.
\end{remark}

\begin{remark}
  \label{rem:normalisation-constant-proj}
  The normalisation constant is particularly simple in the case of projection
  DPPs. Let $\bU\bU^\top$ denote a projection kernel. Then (trivially), $m$ of its
  eigenvalues equal $1$ and the rest are null. By lemma \ref{lem:esp},
  \[ \sum_X \det(\bU\bU^\top)_X = e_m(\bU\bU^\top) = 1 \]
  If as above $\bL = \bU \bM \bM^\top \bU^\top$ with $\bU \in \R^{n \times m}$, then
  by the same reasoning as in the proof of lemma \ref{lem:max-rank-dpp}:
  \[ \sum_X \det(\bL_X) = \det(\bM^\top\bM) \sum_X \det(\bU\bU^\top)_X  = \det(\bM^\top\bM) \]
\end{remark}

\subsubsection{Mixture representation}
\label{sec:mixture-representation}

Determinantal point processes have a well-known representation as a mixture of
projection-DPPs (also sometimes called ``elementary DPPs'' in the literature). See \cite{Barthelme:AsEqFixedSizeDPP}
for details. 
The following mixture representation (due to \cite{Hough:DPPandIndep}) is
fundamental, both for theoretical and computational purposes, since it serves as
the basis for exact sampling of DPPs. There are two variants, one for DPPs and
one for fixed-size DPPs. %For the purposes of this paper, we describe here the mixture representation of L-ensembles only.

\begin{lemma}[Mixture representation of fixed-size L-ensembles \cite{KuleszaTaskar:FixedSizeDPPs}]
  Let $\X \sim \mDPP{m}(\bL)$ be an L-ensemble based on $\bL$, and $\bL = \bU \matr{\Lambda} \bU^\top$ be the spectral decomposition of $\bL$. Then, 
  $\X$ may be obtained from the following mixture process:
  \begin{enumerate}
  \item Sample $m$ indices $\Y \sim \mDPP{m}(\matr{\Lambda})$
  \item Form the projection matrix $\bM = \bU_{:,\Y} (\bU_{:,\Y})^\top$
  \item Sample $\X | \Y \sim  \mDPP{m}(\bM)$
  \end{enumerate}
  Equivalently, the probability mass function of $\X$ can be written as:
  \begin{equation}
    \label{eq:mixture-representation-fixed}
    \Proba(\X = X) =  \frac{\Ind(|X| = m)}{e_m(\matr{\Lambda})} \sum_{Y,|Y|=m}  \det\begin{pmatrix} \bU_{X,Y}   \end{pmatrix}^2 \prod_{i \in Y} \lambda_i
\end{equation}
\end{lemma}

The mixture representation can be understood as (a) first sample which
eigenvectors to use and (b) sample a projection DPP with the selected
eigenvectors.

The counterpart for varying-size L-ensembles looks highly similar. 

\begin{lemma}[Mixture representation of L-ensembles, see e.g. \cite{kulesza2012determinantal}]\label{lem:mixture_dpp}
  Let $\X \sim DPP(\bL)$ and $\bL = \bU \matr{\Lambda} \bU^\top$. Then,
  $\X$ may be obtained from the following mixture process:
  \begin{enumerate}
  \item Sample indices $\Y \sim DPP(\matr{\Lambda})$
  \item Form the projection matrix $\bM = \bU_{:,\Y} (\bU_{:,\Y})^\top$
  \item Sample $\X | \Y \sim  \mDPP{|\Y|}(\bM)$
  \end{enumerate}
  Equivalently, the probability mass function of $\X$ can be written as:
  \begin{equation}
    \label{eq:mixture-representation-varying}
    \Proba(\X = X) = \frac{1}{\det(\bL + \bI)} \sum_{Y}  \det\begin{pmatrix} \bU_{X,Y}   \end{pmatrix}^2 \prod_{i \in Y} \lambda_i.
  \end{equation}
\end{lemma}

The only step that varies is the first one, where we sample from
$DPP(\matr{\Lambda})$ instead of $\mDPP{m}(\matr{\Lambda})$.

\section{Extended L-ensembles: a new representation for DPPs}
\label{sec:ppDPPs}

This introductory section enables us to precisely shed light on what is lacking in the current state-of-the-art regarding explicit joint probability distributions of DPPs.
	
Concerning varying-size DPPs, we observe that i/~DPPs whose marginal kernel do not contain any eigenvalue equal to $1$ are L-ensembles and thus have an explicit joint probability distribution (using Lemma~\ref{lem:LisDPP} and Eq.\eqref{eq:def-dpp_via_L}), ii/~projection DPPs, that is, DPPs whose marginal kernel have eigenvalues equal to exactly $0$ or $1$, are fixed-size L-ensembles, and thus also have an explicit joint probability distribution (using Lemma~\ref{lem:marginal-kernel-proj} and Eq.~\eqref{eq:def-fsdpp}). The current state-of-the-art however does not provide\footnote{A formula due to~\cite{Macchi:CoincidenceApproach} exists in this case but it is unwieldy. See also the discussion around Corollary 1.D.3 in~\cite{Gau20}.} easy-to-use explicit joint distributions for DPPs whose marginal kernels lie in-between those two cases: DPPs whose marginal kernel have some eigenvalues equal to $1$ and others lying strictly between $0$ and $1$. According to Lemma~\ref{lemma:bernoulli}, those are the cases where the size of the DPP is the sum of a deterministic part (the number of eigenvalues equal to 1) and a random part. We will call these particular DPPs \emph{partial projection DPPs} for reasons that will become clear when we study their mixture representation.

Concerning fixed-size DPPs, we observe that all fixed-size L-ensembles are fixed-size DPPs but the contrary is false. In other words, the current state-of-the-art does not provide explicit joint distributions for those fixed-size DPPs that are not fixed-size L-ensembles.

The main contribution of this paper is to provide a unifying framework to write handy, explicit joint distributions for \emph{all} varying and fixed-size DPPs and thus to fill in the holes of the current theory. 
To this end, we introduce \emph{extended L-ensembles}, a novel way of representing the class of DPPs. Before we give their formal definition (Definition~\ref{def:ext_Lens}) we need a few preliminary results.
%With these new tools in hand, one will be able to easily handle joint probability distributions of \emph{any} varying-size or fixed-size DPP.}  %This representation has the advantage of giving explicit expressions for the
%joint probability distribution of \emph{all} varying and fixed-size DPPs.

%In particular, the extended L-ensemble viewpoint will provide easy-to-use, explicit formulas for the joint probability of DPPs in cases where the spectrum of the DPP's marginal kernel contains eigenvalues equal to $1$ (that is, in cases where the DPP at hand is not an L-ensemble) \footnote{A formula due to   \cite{Macchi:CoincidenceApproach} exists in this case but it is unwieldy}.   According to Lemma~\ref{lemma:bernoulli}, those are the cases where the size of the DPP is the sum of a deterministic part (the number of such eigenvalues equal to 1) and a random part. We will call such DPPs \emph{partial projection DPPs} for reasons that will become clear when we study their mixture representation. 

\subsection{Conditionally positive (semi-)definite matrices}
\label{sec:cpd-matrices}

L-ensembles are naturally formed from positive semi-definite matrices, because
$\bL$ being positive semi-definite is a sufficient condition for $\det \bL_\X$
being non-negative. Extended L-ensembles can accomodate a broader set
of matrices called conditionally positive semi-definite (CPD) matrices.
\begin{definition}
  A matrix $\bL \in \R^{n \times n}$ is called conditionally positive
  (semi-)definite with respect to a rank $p\geq 0$ matrix $\bV \in \R^{n \times p}$ if
  $\vect{x}^\top\bL\vect{x} > 0$ (resp., $\vect{x}^\top\bL\vect{x} \geq 0$) for all $\vect{x}$ verifying $\bV^\top \vect{x}=0$.
\end{definition}
\begin{remark}
	Note that we authorize $p=0$ in the definition: in this case, the definition simply boils down to that of positive semi-definite matrices. 
\end{remark}
The set of vectors such that $\bV^\top \vect{x}=0$ is the space orthogonal to the
span of $\bV$, which we note $\orth \bV$. The conditionally positive definite
requirement may be read as a requirement for $\bL$ to be positive definite
within $\orth \bV$. Positive-definite matrices are therefore also conditionally
positive-definite, but matrices with negative eigenvalues may also be
conditionally positive-definite.
\begin{proposition}
  \label{prop:positivity-cond-eigenval}
  Let $\bL$ be conditionally positive (semi-)definite with respect to $\bV \in \R^{n \times p}$, that we suppose full column rank. Let
  $\bQ \in \R^{n \times p}$ designate an orthonormal basis for $\mspan \bV$, so that $\bI - \bQ\bQ^\top$
  is the orthogonal projection on $\orth \bV$. Let $\widetilde{\bL} = (\bI - \bQ\bQ^\top)\bL(\bI -
  \bQ\bQ^\top)$. $\widetilde{\bL}$ is symmetric and thus diagonalisable in $\mathbb{R}$, and all its eigenvalues are non-negative.
\end{proposition}
\begin{proof}
  Follows directly from the definition: $\forall\vect{x}$, $\vect{x}^\top (\bI - \bQ\bQ^\top)\bL(\bI -
  \bQ\bQ^\top)\vect{x} \geq 0$. 
\end{proof}

The above remark will become important when we define extended L-ensembles. The
following example of a conditionally positive definite is classical (but
surprising), and is a special case of a class of conditionally positive definite
kernels studied in \cite{micchelli1986interpolation}. 
\begin{example}[\cite{micchelli1986interpolation}]
  Let
  \[ \bD^{(1)} = [\norm{\vect{x}_i-\vect{x}_j}]_{i,j}\] the distance matrix
  between $n$ points in $\R^d$. Then $-\bD^{(1)}$ is conditionally positive
  definite with respect to the all-ones vector $\vect{1}_n$.
\end{example}
Some extensions of this example can be found in section \ref{sec:ppdpp-examples-cpdef}. 
\subsection{Nonnegative Pairs}
\label{sec:ppDPP-def}

The central object when defining extended L-ensembles is what we call a Nonnegative Pair (NNP for short). 

\begin{definition}
	\label{def:nnp}
  A Nonnegative Pair, noted $\ELE{\bL}{\bV}$ is a pair $\bL \in \R^{n \times n}$, $\bV \in \R^{n
    \times p}$, $0\leq p\leq n$, such that $\bL$ is symmetric and conditionally positive semi-definite with respect to
  $\bV$, and $\bV$ has full column rank. Wherever a NNP $\ELE{\bL}{\bV}$ appears below, we consistently use the following notation: 
  \begin{itemize}
  		\item $\bQ \in \R^{n \times p}$ is an
  		orthonormal basis of  $\mspan \bV$, such that $\bI - \bQ\bQ^\top$
  		is a projector on $\orth \bV$
  		\item $\widetilde{\bL} = (\bI - \bQ\bQ^\top)\bL(\bI -
  		\bQ\bQ^\top)\in \R^{n \times n}$. From Proposition ~\ref{prop:positivity-cond-eigenval}, we know that all its eigenvalues are non-negative. We will denote by $q$ the rank of $\widetilde{\bL}$. Note that $q \le n-p$ as the $p$ columns of $\bQ$ are trivially eigenvectors of $\widetilde{\bL}$ associated to $0$. We write
  		\[
  		\widetilde{\bL} = \tbU\tbLam\tbU^\top
  		\] 
  		its truncated spectral decomposition; where $\tbLam=\text{diag}(\widetilde{\lambda}_1,\ldots,\widetilde{\lambda}_q)\in \mathbb{R}^{q\times q}$ and $\tbU \in \mathbb{R}^{n\times q}$ are the diagonal matrix of nonzero eigenvalues and the matrix of the corresponding eigenvectors of $\widetilde{\bL}$, respectively.
  \end{itemize}
\end{definition}
\begin{remark}
		Again, note that we authorize $p=0$ in the definition: in this case, $\tbL$ boils down to $\bL$.
\end{remark}
Let us now formulate the following lemma, useful for the next section. 
\begin{lemma}\label{lem:equivalence_ppdpp_ensemble}
	Let $\ELE{\bL}{\bV}$ be a NNP. Then, for any subset $X \subseteq \{1, \ldots, n\}$:
	\begin{align*}
	(-1)^{p}\det \begin{pmatrix}
	{\bL}_{X} & \bV_{X,:} \\
	(\bV_{X,:})^\top & \matr{0} 
	\end{pmatrix} = (-1)^{p} \det \begin{pmatrix}
	{\tbL}_{X} & \bV_{X,:} \\
	(\bV_{X,:})^\top & \matr{0} 
	\end{pmatrix}  \ge 0.
	\end{align*}
\end{lemma}
\begin{proof} 
	Let us write $m=|X|$ the size of $X$. 
	The case $\rank \bV_{X,:} < p$ is trivial as both sides of the equality are zero. Next, assume that $\bV_{X,:}\in\mathbb{R}^{m\times p}$ is full column rank.
	If $m = p$, then $\bV_{X,:}$ is square and both sides are equal to $(\det \bV_{X,:})^2$ by lemma~\ref{lem:block-det}. 
	Now consider the case $m>p$.
	Let $\bQ$ be as in Definition~\ref{def:nnp}, so that $\bV  = \bQ \matr{R}$ (with $\matr{R}$ nonsingular).  
	Let $\bB(X) \in \R^{m \times (m-p)}$ be a basis of $\orth(\bV_{X,:}) = \orth(\bQ_{X,:}) $. 
	Then, using lemma~\ref{lem:det-saddlepoint}, we have that
	\begin{align*}
	(-1)^p  \det \begin{pmatrix}
	\bL_{X} & \bV_{X,:} \\
	(\bV_{X,:})^\top & \matr{0} 
	\end{pmatrix}  &= \det ( (\bV_{X,:})^{\top}  \bV_{X,:})  \det ((\bB^\top(X) {\bL}_{X}\bB(X))\\
	&= \det ( (\bV_{X,:})^{\top}  \bV_{X,:})  \det ((\bB^\top(X) \tbL_{X}\bB(X) )
	\\
	&=  (-1)^p\det
	\begin{pmatrix}
	\widetilde{\bL}_{X} & \bV_{X,:} \\
	(\bV_{X,:})^\top & \matr{0} 
	\end{pmatrix}, 
	\end{align*}
	where the last but one equality is from   $\tbL=(\bI - \bQ\bQ^\top)\bL(\bI - \bQ\bQ^\top)$ and the fact that  $\bB^{\T}(X)\bQ_{X,:} = 0$ and hence $(\bI - \bQ\bQ^\top)_X \bB(X) = \bB(X)$. 
	Finally, $\det (\bB^\top(X) \tbL_{X}\bB(X) ) \ge 0$ due to positive semidefiniteness of $\tbL$, which completes the proof.
\end{proof}

\subsection{DPPs via extended L-ensembles}

\begin{definition}[Extended L-ensemble] 
	\label{def:ext_Lens}
	Let $\ELE{\bL}{\bV}$ be any NNP. An extended L-ensemble $\X$ based on $\ELE{\bL}{\bV}$ is a point process verifying:
	\begin{align}
	\label{eq:proba_mass}
	\forall X\subseteq\Omega,\qquad \Proba(\X=X) \propto (-1)^p \det
	\begin{pmatrix}
	\bL_{X} & \bV_{X,:} \\
	(\bV_{X,:})^\top & \matr{0}
	\end{pmatrix}.
	\end{align} 
\end{definition}
	\begin{remark}
	We stress that an extended L-ensemble reduces to an L-ensemble only in the case $p=0$. If $p\geq 1$, an extended L-ensemble is not an L-ensemble, since the probability mass function of $\X$ is not expressed as a principal minor of a larger matrix\footnote{another way to see that
		 is by looking at the distribution of the size of the sampled set: by definition~\ref{def:dpp_via_L}, an L-ensemble has a non-null probability of sampling the empty set; whereas Corollary~\ref{cor:distrib_size_ELE} states that the empty set has a null probability of being sampled as soon as $p\geq 1$}. Also, 
	the right-hand  side in eq. \eqref{eq:proba_mass} is non-negative by
	Lemma~\ref{lem:equivalence_ppdpp_ensemble}, and thus defines a valid
	probability distribution. The normalisation constant is tractable and given
	later (see section \ref{sec:ppdpp-marginal-kernel}). 
	On a more minor note, the factor $(-1)^p$ arises because of the peculiar properties of saddle-point matrices, see Lemma~\ref{lem:det-saddlepoint}. 
\end{remark}

Importantly, the class of extended L-ensembles is identical to the class of DPPs, as the two following theorems demonstrate. 
	\begin{theorem}
		\label{thm:L-ens_to_K}
		Let $\ELE{\bL}{\bV}$ be any NNP, and $\X$ be an extended L-ensemble based on $\ELE{\bL}{\bV}$. 
		Then, $\X$ is a DPP with marginal kernel
		\begin{equation}
		\label{eq:dpDPPmarginalkernel}
		\bK = \bQ \bQ^\top + \tbL(\bI+ \tbL)^{-1}.\end{equation}
	\end{theorem}
	\begin{proof}
		See Appendix~\ref{sec:thm:K_to_L-ens_and_back}.

\end{proof}
The converse is also true: any DPP (not only L-ensembles) is an extended L-ensemble:

	\begin{theorem}
		\label{thm:K_to_L-ens}
		Let $\bm{0}\preceq\bK\preceq \bI$ be any marginal kernel and $\X$ its associated DPP. Denote by $\bV\in\mathbb{R}^{n\times p}$ the matrix concatenating the $p\geq 0$ orthonormal eigenvectors of $\bK$ associated to eigenvalue $1$ and $\bL=\bK\left(\bI-\bK\right)^{\dagger}$ with $\dagger$ representing the Moore-Penrose pseudo-inverse. Then, $\X$ is an extended L-ensemble based on the NNP $\ELE{\bL}{\bV}$.
	\end{theorem}
	\begin{proof}
		See Appendix~\ref{sec:thm:K_to_L-ens_and_back}.
	\end{proof}
Recall that, as per definition~\ref{def:fsDPP}, a fixed-size DPP is simply a DPP conditioned on size. As a consequence of the equivalence between extended L-ensembles and DPPs, one obtains the following explicit expression of the probability mass function of any fixed-size DPP:
	\begin{corollary}
		Let $\bm{0}\preceq\bK\preceq \bI$ be any marginal kernel and $\X$ its associated fixed-size DPP of size $m$. Let $\ELE{\bL}{\bV}$ be the NNP as defined in theorem~\ref{thm:K_to_L-ens}. Then 
		\begin{align}
		\label{eq:proba_mass_fs}
		\forall X\subseteq\Omega,\qquad \Proba(\X=X) \propto (-1)^p \det
		\begin{pmatrix}
		\bL_{X} & \bV_{X,:} \\
		(\bV_{X,:})^\top & \matr{0}
		\end{pmatrix}\Ind(|X|=m).
		\end{align}
	\end{corollary}
	\begin{remark}
		Fixed-size DPPs of size $m$ with marginal kernel $\bK$ cannot be defined for
    $m$ smaller than the multiplicity of $1$ in the spectrum of $\bK$ (by lemma~\ref{lemma:bernoulli}). 
    Consequently, from the extended L-ensemble viewpoint, $m$ should always be
    larger than or equal to $p$.
	\end{remark}

\subsubsection{Partial projection DPPs}
%The previous section made clear that 
%	\begin{itemize}
%		\item any DPP in the class of DPPs may be defined equivalently either via a marginal kernel  $\bm{0}\preceq\bK\preceq\bI$ from the marginal point of view, or via a NNP $\ELE{\bL}{\bV}$ from the point of view of the explicit probability mass function.
%		\item the class of fixed-size DPPs, being in all generality defined as DPPs conditioned on size, are in fact best described with extended L-ensembles. Their probability mass function are given by Eq.~\eqref{eq:proba_mass_fs}. Apart from the special case where $m=p$ that implies a projection DPP~\footnote{If $m=p$, 
%			$\bV_{X,:}$ is square in Eq.~\eqref{eq:proba_mass_fs} and by Lemma~\ref{lem:det-saddlepoint}, $\Proba(\X = X) \propto \det(\bV_{X,:})^2$, which is the probability mass function of a projection DPP (see lemma \ref{lem:max-rank-dpp}).}, fixed-size DPPs do not have marginal kernels.
%	\end{itemize}
	We differentiate DPPs (both varying-size and fixed-size) defined by NNPs $\ELE{\bL}{\bV}$ for which
	\begin{itemize}
		\item $p=0$:  Eq.~\ref{eq:proba_mass} (resp. Eq.~\ref{eq:proba_mass_fs}) boils down to Eq.~\ref{eq:def-dpp_via_L} (resp. Eq.~\ref{eq:def-dpp}): we recover the L-ensembles $\X\sim DPP(\bL)$ (resp. fixed-size L-ensembles $\X\sim|DPP|_m(\bL)$).
		\item $p\geq1$: in this case, the associated DPPs are not L-ensembles; and we will call them \emph{partial-projection DPPs} (pp-DPPs) for reasons that will become clear in section~\ref{sec:mixture_rep_ppDPP}. We will denote them $\X\sim \ppDPP \ELE{\bL}{\bV}$ and $\X\sim \mppDPP{m}\ELE{\bL}{\bV}$ for the varying-size and the fixed-size cases respectively.
\end{itemize}

This ends the first part of our work: the formalism of extended L-ensembles enables to fill in the holes of the current theory of DPPs by providing explicit joint probability distributions for all varying and fixed-size DPPs. We now move on to studying this novel object.

\section{Extended L-ensembles: mixture representation, main properties, sampling}
\label{sec:general_properties}
In this section, we start by showing an extension of the Cauchy-Binet formula to saddle-point matrices (section~\ref{subsec:CB}) that will prove useful to give the mixture representation of extended L-ensembles (section~\ref{sec:mixture_rep_ppDPP}).  In section~\ref{sec:ppdpp-properties}, we list a few basic properties of extended L-ensembles (such as normalisation and complements). Section~\ref{sec:ppdpp-examples} gives a few examples for illustration purposes and section~\ref{sec:ppdpp-sampling} discusses sampling strategies of $\X\sim\ppDPP\ELE{\bL}{\bV}$. 

\subsection{A generalisation of the Cauchy-Binet Formula}
\label{subsec:CB}
The cornerstone of the mixture representation of L-ensembles, discussed in
Section~\ref{sec:mixture-representation}, is in fact the Cauchy-Binet formula,
recalled in Lemma~\ref{lem:cauchy-binet} (see for
instance~\cite{Hough:DPPandIndep, kulesza2012determinantal}). In order to
provide a similar spectral understanding of extended L-ensembles, we need the
following generalisation of  Cauchy-Binet. 

\begin{theorem}[Generalisation of Cauchy-Binet]
	\label{thm:equivalence-extended-spectral}
	Let $\ELE{\bL}{\bV}$ be a NNP, and $\bQ$, $\tbU$, $\tbLam =
  \diag(\widetilde{\lambda}_1 \ldots \widetilde{\lambda}_{q})$ and $q$ be as in Definition~\ref{def:nnp}. Then for any subset $X \subseteq \{1, \ldots, n\}$ of size $|X| = m$, $p\leq m \leq p+q$, 
	it holds that
	\begin{equation}
	\label{eq:equivalence-spectral-extended-pDPP}
	(-1)^p\det
	\begin{pmatrix}
	\bL_{X} & \bV_{X,:} \\
	(\bV_{X,:})^\top & \matr{0} 
	\end{pmatrix} = \det(\bV^\top\bV)\sum_{Y,|Y|=m-p} \det \left( 
	\begin{bmatrix}
	\bQ_{X,:} & \tbU_{X,Y}
	\end{bmatrix}  \right)^2
	\prod_{i \in Y} \widetilde{\lambda}_i
	\end{equation}
\end{theorem}
\begin{proof}
	First of all, writing the $(\bQ, \bR)$ decomposition of $\bV$ as $\bV = \bQ\bR$ one has:
	\[
	\det
	\begin{pmatrix}
	\bL_{X} & \bV_{X,:} \\
	(\bV_{X,:})^\top & \matr{0} 
	\end{pmatrix} 
	= (\det(\matr{R}))^2
	\det
	\begin{pmatrix}
	\bL_{X} & \bQ_{X,:} \\
	(\bQ_{X,:})^\top & \matr{0} 
	\end{pmatrix}.     
	\]
	Noting that $\det(\bV^\top\bV) = (\det(\matr{R}))^2$, to prove Eq.~\eqref{eq:equivalence-spectral-extended-pDPP}   it is sufficient to show that:
	\begin{equation}
	\label{eq:equivalence-spectral-extended-pDPP_Q}
	(-1)^p\det
	\begin{pmatrix}
	\bL_{X} & \bQ_{X,:} \\
	(\bQ_{X,:})^\top & \matr{0} 
	\end{pmatrix} = \sum_{Y,|Y|=m-p} \det \left( 
	\begin{bmatrix}
	\bQ_{X,:} & \tbU_{X,Y}
	\end{bmatrix}  \right)^2
	\prod_{i \in Y} \widetilde{\lambda}_i.
	\end{equation}
	Now, the case $\rank \bQ_{X,:} < p$ is trivial as both sides in \eqref{eq:equivalence-spectral-extended-pDPP_Q} are zero. Next, we assume that $\bQ_{X,:}$ is full rank.
	Using first lemma~\ref{lem:equivalence_ppdpp_ensemble} and then lemma~\ref{lem:det-coef-polynomial}, one has:
	\[
	(-1)^p\det    \begin{pmatrix}  {\bL}_{X} & {\bQ}_{X,:} \\  ({\bQ}_{X,:})^\top & {0}  \end{pmatrix}  =
	(-1)^p\det  \begin{pmatrix}     \widetilde{\bL}_{X} & {\bQ}_{X,:} \\      ({\bQ}_{X,:})^\top & {0}    \end{pmatrix}
	= [t^p ] \det(\widetilde{\bL}_X + t {\bQ}_{X,:}({\bQ}_{X,:})^\top).
	\]
	Using the fact that $\widetilde{\bL} = {\tbU}\tbLam{\tbU}^{\top}$, the right hand side may be re-written:
	\begin{align*}
    [t^p]\det(\widetilde{\bL}_X + t {\bQ}_{X,:}({\bQ}_{X,:})^\top) & =  [t^p]\det\left([{\bQ}_{X,:} {\tilde \bU}_{X,:} ] \begin{pmatrix}t \bI_p & \matr{0} \\ \matr{0} &  \widetilde{\matr{\Lambda}}  \end{pmatrix}[{\bQ}_{X,:} {\tilde \bU}_{X,:}]^\top\right)\\ 
	&  = 
	\sum\limits_{|Y| = m-p} (\det([{\bQ}_{X,:} {\tilde\bU}_{X,Y} ]))^2 \det(\widetilde{\matr{\Lambda}}_{Y}),
	\end{align*}
	where the last equality follows from the  Cauchy-Binet lemma.
\end{proof}

\subsection{Mixture representation}
\label{sec:mixture_rep_ppDPP}

In the mixture representation of L-ensembles (see Sec.~\ref{sec:mixture-representation}), one first samples a set of orthonormal vectors, forms a
projective kernel from these eigenvectors, and then samples a projection DPP
from that kernel. In that sense, a projection DPP is the trivial mixture in
which the same set of eigenvectors is always sampled. 
In this section, we will see that in partial projection DPPs, a subset of orthogonal vectors is included deterministically (coming from $\bV$), and the rest are subject to sampling, from the
part of $\bL$ orthogonal to $\bV$, hence the name \emph{partial projection}. 

In fact, examining Eq.~\eqref{eq:equivalence-spectral-extended-pDPP}, the kinship with the mixture representation of fixed-size L-ensembles should be clear 
upon comparison with equation \eqref{eq:mixture-representation-fixed}. The left-hand side of Eq.~\eqref{eq:equivalence-spectral-extended-pDPP} is the
probability mass function, and on the right-hand side we recognise a sum (over $Y$) of
probability mass functions for projection DPPs ($\det \left( 
\begin{bmatrix}    \bQ_{X,:} & \tbU_{X,Y}  \end{bmatrix}  \right)^2$) indexed by $Y$, weighted by a
product of eigenvalues ($\prod_{i \in Y} \widetilde{\lambda}_i$). This lets us represent the partial-projection DPP as a
probabilistic mixture. Contrary to fixed-size L-ensembles, some eigenvectors
appear with probability 1: the ones that originate from $\bV$ (represented by
$\bQ_{\X,:}$ in Eq.~\eqref{eq:equivalence-spectral-extended-pDPP}). The rest are
picked randomly according to the law given by the product $\Proba(\Y = Y) \propto \prod_{i \in Y} \tilde{\lambda}_i$.

Seen as a statement about probabilistic mixtures,  theorem~\ref{thm:equivalence-extended-spectral} provides a recipe
for sampling from $\X \sim \mppDPP{m} \ELE{\bL}{\bV} $. We summarize this recipe in the following statement: 

\begin{corollary}
  \label{cor:mixture-rep-ppdpp}
  Let $\ELE{\bL}{\bV}$ be a NNP, and $\bQ$, $\tbU$, $\tbLam$ and $q$ be as in Definition~\ref{def:nnp}. 
  Let  $\X \sim\mppDPP{m} \ELE{\bL}{\bV}$ with $p \leq  m \leq  p+q$. Then,
  $\X$ may be obtained from the following mixture process:
  \begin{enumerate}
  \item Sample $m-p$ indices $\Y \sim \mDPP{m-p}(\tbLam)$
  \item Form the projection matrix $\bM = \bQ\bQ^\top + \tbU_{:,\Y}
    (\tbU_{:,\Y})^\top$ (recall that $\bQ$ and $\tbU$ are orthogonal)
  \item Sample $\X | \Y \sim  \mDPP{m}(\bM)$
  \end{enumerate}
\end{corollary}
Note that at step 1 we only sample from the \emph{optional} part, since the
eigenvectors from $\bV$ need to be included anyway. The total number of
eigenvectors to include is $m$, so $m-p$ need to be sampled randomly. 

Using theorem \ref{thm:equivalence-extended-spectral}, as in the fixed-size
case, we arrive easily at the following mixture characterisation for the varying-size case: 

\begin{corollary}
  \label{cor:mixture-rep-ppDPP-varying}
  Let $\ELE{\bL}{\bV}$ be a NNP, and $\bQ$, $\tbU$ and $\tbLam$ be as in Definition~\ref{def:nnp}. 
  Let $\X \sim \ppDPP \ELE{\bL}{\bV}$. Then,
  $\X$ may be obtained from the following mixture process:
  \begin{enumerate}
  \item Sample indices $\Y \sim DPP(\tbLam)$
  \item Form the projection matrix $\bM = \bQ\bQ^\top + \tbU_{:,\Y} (\tbU_{:,\Y})^\top$
  \item Sample $\X | \Y \sim  \mDPP{p+|\Y|}(\bM)$
  \end{enumerate}
\end{corollary}

The only difference from the fixed-size case is in step 1. Again, we include all eigenvectors from $\bV$
(they make up the $\bQ\bQ^\top$ part of the projection matrix $\bM$),
then the remaining ones are sampled from $\Y \sim
DPP(\tbLam)$, which is equivalent to including the eigenvector $\tilde{\vect{u}}_i$
with probability $\frac{\tilde{\lambda}_i}{1+\tilde{\lambda}_i}$.

\subsection{Properties}
\label{sec:ppdpp-properties}

\subsubsection{Normalisation}
\label{sec:ppdpp-marginal-kernel}

Using theorem \ref{thm:equivalence-extended-spectral}, the normalisation
constant is tractable both in the fixed-size and varying-size cases, as shown by the following corollary (see also \cite[Lemma 3.11]{BarthelmeUsevich:KernelsFlatLimit} for an alternative formulation).

\begin{corollary}
\label{cor:normalisation-ppDPP}
  Let $\ELE{\bL}{\bV}$ be a NNP, and $\widetilde{\bL}$ and $q$ as in Definition~\ref{def:nnp}. For $m$ such that $p\leq m \le n$, one has:
  \begin{equation}
    \label{eq:marginalisation-ppDPP}
    (-1)^p \sum_{|X|=m} \det
    \begin{pmatrix}
      \bL_{X} & \bV_{X,:} \\
      (\bV_{X,:})^\top & \matr{0} 
    \end{pmatrix} = e_{m-p}(\tbL)\det(\bV^\top\bV)
  \end{equation}
  and
  \begin{equation}
  (-1)^p \sum_{X} \det
    \begin{pmatrix}
      \bL_{X} & \bV_{X,:} \\
      (\bV_{X,:})^\top & \matr{0} 
    \end{pmatrix} = \det(\bI + \tbL)\det(\bV^\top\bV)
    \label{eq:marginalisation-ppDPP-varying}
\end{equation}

\end{corollary}
\begin{proof}
If $m > p+q$, then the right-hand side is zero, as well as the left-hand side (by lemma~\ref{lem:det-saddlepoint}).  
 In the case $m \le p+q$,  from  theorem \ref{thm:equivalence-extended-spectral} we have:
  \begin{align*}
    (-1)^p \sum_{|X|=m} \det
    \begin{pmatrix}
      \bL_{X} & \bV_{X,:} \\
      (\bV_{X,:})^\top & \matr{0} 
    \end{pmatrix} &= 
                     \det(\bV^\top\bV) \sum_{|X|=m} \sum_{Y,|Y|=m-p} \det \left( 
    \begin{bmatrix}
      \bQ_{X,:} & \tbU_{X,Y}
    \end{bmatrix}  \right)^2
                   \prod_{i \in Y} \tilde{\lambda}_i \\
              &=  \det(\bV^\top\bV) \sum_{Y,|Y|=m-p} \quad \prod_{i \in Y} \tilde{\lambda}_i \\
              &= e_{m-p}(\tilde{\bL})\det(\bV^\top\bV),
  \end{align*}
  where the sum over $X$ is just the normalisation constant of a projection DPP
  (see remark \ref{rem:normalisation-constant-proj}). The proof for varying size
  is similar, using:
  $\sum_{Y} \prod_{i \in Y} \tilde{\lambda}_i = \prod_{i=1}^q (1 + \tilde{\lambda}_i).$
\end{proof}

Using these results, we easily obtain the distribution of the size of $|\X|$
for $\X\sim \ppDPP \ELE{\bL}{\bV}$. One may check that equivalent results
are obtained either using the mixture representation (see corollary
\ref{cor:mixture-rep-ppDPP-varying}) or the associated marginal kernel (via Eq.~\ref{eq:dpDPPmarginalkernel} and lemma~\ref{lemma:bernoulli}). 
\begin{corollary}\label{cor:distrib_size_ELE}
  Let $\X \sim \ppDPP \ELE{\bL}{\bV}$. Then
  \begin{equation}
    \label{eq:prob_size_ppdpp}
    \Proba(|\X|=m) =
    \begin{cases}
      0, & \mbox{\ if\ } m < p, \\
      \frac{e_{m-p}(\tilde{\bL})}{\det(\tbL + \bI)}, &\mbox{\ otherwise}.
    \end{cases}
  \end{equation}
\end{corollary}

\subsubsection{Complements of DPPs}
\label{sec:complements_dpps}

A known (see e.g., \cite{kulesza2012determinantal}, section 2.3) result about DPPs is
that the complement of a DPP in $\Omega$ is also a DPP, i.e., if $\X$ is a DPP,
$\X^c = \Omega \setminus \X$ is also a DPP. We shall give a short proof and some
extensions. 

\begin{theorem}
  \label{thm:complement-DPP}
  Let $\X$ be a DPP with marginal kernel $\bK$. Then the complement of
  $\X$, noted $\X^c$, is also a DPP, and its marginal kernel is $\bI - \bK$. 
\end{theorem}
\begin{proof}
  We first prove this for projection DPPs. Let $\cA \sim \mDPP{m}(\bU\bU^\top)$ for
  orthogonal $\bU$ of rank $m$. Then
  \[ \Proba(\cA^c=A) = \Proba(\cA = A^c) \propto \det( \bU_{A^c,:})^2.\] Note that for the
  probability to be non null we need $A$ to be of size $ n-m $.
  
  Let $\bV \in \R^{n \times (n-m)}$ so that $ \bI = \bU\bU^\top + \bV\bV^\top$. $ \bM
  = 
  \begin{pmatrix}
    \bU & \bV
  \end{pmatrix}$ is an orthogonal basis for $\R^n$ which we may partition as
  $\begin{pmatrix}
    \bU_{A^c,:} & \bV_{A^c,:} \\
    \bU_{A,:} & \bV_{A,:}
  \end{pmatrix}$.
  By lemma \ref{lem:block-det}
  \[ \det \bM  =  \det \bU_{A^c,:} \det \left( \bV_{A,:} -
      \bU_{A,:} (\bU_{A^c,:})^{-1}  \bV_{A^c,:} \right). \]
  This gives
  \[ \Proba(\cA^c = A) \propto \det \left( (\bV_{A,:} - 
      \bU_{A,:} (\bU_{A^c,:})^{-1}  \bV_{A^c,:}) \right)^{-2}.\]
  By the inversion formula for block matrices this is equal to the lower-right
  block in $\bM^{-1} =\bM^\top$, and so:
  \[ \Proba(\cA^c = A) \propto \det \left( \bV_{A,:} \right)^2\] where we recognise a
  projection DPP ($\cA^c \sim \mDPP{n-m}(\bV\bV^\top)$, as claimed). 
  We now use the mixture property to show the general case.
  In the general case,
  \[ \Proba(\X = X) = \sum_{\Y} \Proba(\Y) \det (\bU_{X,\Y})^2\]
  so that:
  \[ \Proba(\X^c = A) = \sum_{\Y} \Proba(\Y) \det (\bU_{A^c,\Y})^2
    =  \sum_{\Y^c} \Proba(\Y^c) \det (\bV_{A,\Y^c})^2 \]
  Since each eigenvector is picked independently in $\Proba(\Y)$ with probability
  $\pi_i$, picking each eigenvector independently with probability $1-\pi_i$
  produces a draw from $\Proba(\Y^c)$. $\Proba(\X^c = A)$ is therefore a DPP, and its kernel is $\bI-\bK$.
\end{proof}

Applying the theorem to L-ensembles we obtain: 
\begin{corollary}
  Let $\X \sim DPP(\bL)$, with $\bL$ a rank $r$ matrix and $r \leq n$. Then
  $\X^c \sim DPP(\bL^\dag, \bV)$ with $\bV$ a basis for $\orth \bL$. In
  particular, if $r = n$
  ($\bL$ is full rank), we have $\X^c \sim DPP(\bL^{-1})$.
\end{corollary}

For extended L-ensembles this generalises to:
\begin{corollary}
  Let $\X \sim DPP \ELE{\bL}{\bV} $, and let $\bZ$ be a basis for $\orth \tbL
  \setminus \mspan \bV $. Then 
  $\X^c \sim DPP(\tbL^\dag, \bZ)$. 
\end{corollary}

The following fixed-size variant is new: it states that the complement
of a fixed-size DPP is also a fixed-size DPP
\begin{proposition}
  Let $\X \sim \mDPP{m} \ELE{\bL}{\bV} $, and let $\bZ$ be a basis for $\orth \tbL  \setminus \mspan \bV $. Then 
  $\X^c \sim \mDPP{n-m}(\tbL^\dag, \bZ)$.
\end{proposition}
\begin{proof}
  Proof sketch: repeat the proof of th. \ref{thm:complement-DPP} up to the mixture
  representation, where we note that since $p(\Y = Y) \propto \prod_{i \in Y} \lambda_i$,
  $p(\Y^c) \propto \prod_{j \in Y^C} \frac{1}{\lambda_j}$ which is again a
  diagonal fixed-size DPP. 
\end{proof}

\subsubsection{Partial Invariance}
\label{sec:ppdpp-invariance}

We parametrise partial-projection DPPs using a pair of matrices (the NNP $\ELE{\bL}{\bV}$), but this is an over-parameterisation since all
that matters is the linear space spanned by $\bV$, as the following makes clear: 
\begin{remark}
	\label{rem:first_inv}
  Consider a NNP $\ELE{\bL}{\bV}$. Let $\bV'=\bV\bR$ with $\bR \in \R^{p
    \times p}$ invertible. We have $\mspan \bV' = \mspan \bV$. Then $\X
  \sim \ppDPP\ELE{\bL}{\bV}$ and  $\X
  \sim \ppDPP\ELE{\bL}{\bV'}$  define the same point process.
  This also holds for $\X \sim \mppDPP{m}\ELE{\bL}{\bV}$ for any $m \geq p$.
\end{remark}
\begin{proof}
  This is clear from theorem \ref{thm:equivalence-extended-spectral} or the
  mixture representation of the partial-projection DPP. Nothing on the
  right-hand side of equation \eqref{eq:equivalence-spectral-extended-pDPP}
  is affected by replacing $\bV$ with a matrix with identical span. In
  particular, the distribution is invariant to rescaling of $\bV$ by any
  non-zero scalar. 
\end{proof}

Notice that this generalises a property of projection DPPs given in the
introduction (section \ref{sec:diag-and-proj}), which is that $\X \sim
\mDPP{m}(\bL)$ and $\X \sim \mDPP{m}(\bL')$ are the same if $\bL$ and $\bL'$
have the same column span and rank $m$. 

Another source of invariance in partial projection DPPs lies in $\bL$: we can
modify $\bL$ along the subspace spanned by $\bV$ without changing the distribution.
\begin{remark}
	\label{rem:second_inv}
  Consider a NNP $\ELE{\bL}{\bV}$. Let $\bL' = \bL + \bV\bm{X}^\top + \bm{X}\bV^\top$ for any matrix $\bm{X} 
  \in \R^{n \times p}$. Then $\X \sim \ppDPP\ELE{\bL}{\bV}$ and $\X' \sim
  \ppDPP\ELE{\bL'}{\bV}$ have the same distribution.
\end{remark}
\begin{proof}
Indeed, by Definition~\ref{def:ext_Lens}, we have $\widetilde{\bL' } = \tbL$.
Therefore, by lemma~\ref{lem:equivalence_ppdpp_ensemble}, the DPPs defined by  ${\bL' }$ and $\bL$ coincide. 
\end{proof}

\subsection{Examples}
\label{sec:ppdpp-examples}

We give here a few examples of partial projection DPPs and their NNPs.

\subsubsection{Partial projection DPPs as conditional distributions}

A simple example of a partial projection DPP arises when the columns of the matrix $\bV$ come from a canonical basis (i.e., each column of $\bV$ is a standard unit vector).
In this case, partial projection DPPs can be interpreted as a particular conditioning of a DPP.
For simplicity, assume that $\matr{V}  =  \begin{bmatrix} \matr{I}_p & 0 \end{bmatrix}^{\top}$,
so that the projected $\matr{L}$ matrix  becomes
\[
\tbL =  \begin{bmatrix} 0& 0  \\ 0 & \bL_{\{\vect{x}_{p+1}, \ldots,\vect{x}_n\}}\end{bmatrix}.
\]

In this case, the mixture representation for partial projection DPPs (resp. fixed-sized  partial projection DPPs) implies that:
\begin{itemize}
\item all the points  $\vect{x}_1, \ldots,\vect{x}_p$ are always sampled; 
\item  the remaining points are sampled according to the L-ensemble (resp. fixed-size L-ensemble) based on   $\matr{L}_{\{\vect{x}_{p+1}, \ldots,\vect{x}_n\}}$.
\end{itemize}
For example, in the  varying-size case $\X \sim \ppDPP \ELE{\bL}{\bV}$, the probability of sampling the remaining points is
\begin{equation}\label{eq:pp-dpp-conditional}
\Proba(\X \cap \{\vect{x}_{p+1}, \ldots,\vect{x}_n\} = X') \propto \det  \bL_{X'},
\end{equation}
which is linked  to a certain conditional distribution of the ordinary L-ensemble based on $\bL$
(see \cite[\S 2.4.3]{kulesza2012determinantal} for more details).

\subsubsection{Roots of trees in uniform spanning random forests are partial projection DPPs}
\label{sec:roots-forests}

\begin{figure}[pt]
\begin{center}
\begin{picture}(0,0)\includegraphics{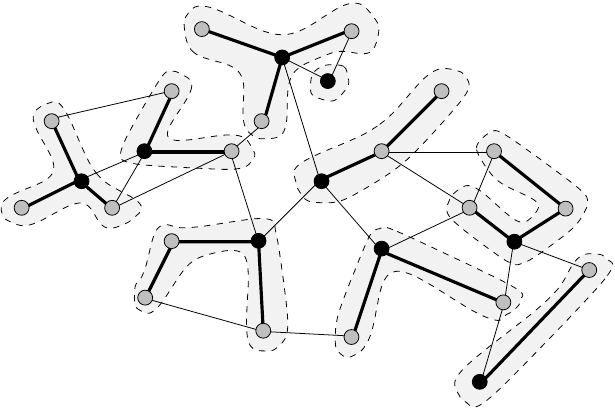}\end{picture}\setlength{\unitlength}{3315sp}\begingroup\makeatletter\ifx\SetFigFont\undefined \gdef\SetFigFont#1#2#3#4#5{\reset@font\fontsize{#1}{#2pt}\fontfamily{#3}\fontseries{#4}\fontshape{#5}\selectfont}\fi\endgroup \begin{picture}(5850,3872)(2937,-4753)
\end{picture} \end{center}
\caption{Roots of uniform random forests over a graph are distributed according to a partial projection DPP. Vertices or nodes are depicted in gray; edges as thin lines. A random forest is depicted : its trees are surrounded by light gray zones; edges of the trees are thicks black lines; roots are the black nodes. The forest is spanning the graph as each nodes of the graph appears once in a tree of the forest.}
\label{fig:URForest}
\end{figure}

It is known ({\it e.g. } \cite{avena2013some}) that the roots of the trees in a uniform random spanning forest over a graph with $n$ nodes and Laplacian $\vect{{\mathcal{L}}}$ are distributed according to a DPP with  marginal kernel  $\vect{K} = q (q \vect{I} +\vect{{\mathcal{L}}})^{-1} $ for some real parameter $q>0$. 
Figure \ref{fig:URForest} illustrates what a spanning forest over a graph is. 
Let us denote as $\lambda_1\geq \ldots \geq \lambda_n = 0$ the eigenvalues of the Laplacian, and $\{ \vect{u}_i \}_i$ the associated set of orthonormal  eigenvectors. It is well known that $\lambda_n=0$ for any graph:  $\vect{K}$ has  thus at least one eigenvalue equal to $1$ and, as such, the associated DPP is not an L-ensemble. It can however be described by an extended L-ensemble:
\begin{proposition}
The set of roots in a uniform random spanning forest over a connected graph with Laplacian $\vect{{\mathcal{L}}}$ is distributed according to a partial projection DPP with NNP $( q \vect{{\mathcal{L}}}^{\dagger} ; \vect{1} )$, where $\dagger$ stands for the Moore-Penrose inverse.
\end{proposition}

\begin{proof}
Applying theorem~\ref{thm:K_to_L-ens}, a DPP with marginal kernel $\bK$
  can be described by an extended L-ensemble based on the NNP $(\bL, \bV)$ with $\bV$ and $\bL$ verifying:
	\begin{itemize}
		\item the matrix $\bV$ concatenates all eigenvectors of $\bK$ associated to eigenvalue 1: in a connected graph, there is only one such eigenvalue and it is associated to eigenvector  $\vect{u}_n=n^{-1/2} \vect{1}$
		\item the matrix $\bL$ is equal to $\bK(\bI-\bK)^\dagger$, which is equal to $q\vect{{\mathcal{L}}}^{\dagger}$
	\end{itemize}
\end{proof}

\begin{remark}
  This example also provides a nice illustration for the properties of
  complements of DPPs (section \ref{sec:complements_dpps}). Since $\vect{{\mathcal{L}}}$ is a positive-definite matrix, we may define  $\mathcal{C} \sim
  DPP(\frac{1}{q}\vect{{\mathcal{L}}})$. The complement of $\mathcal{C}$ is a DPP $\mathcal{C}^c
  \sim DPP \ELE{q\vect{{\mathcal{L}}}^\dag}{\mathbf{1}}$, which  from the result above 
  corresponds to the roots process. 
  $\mathcal{C}$ therefore samples every node except the roots of a random forest
  on the graph. 
\end{remark}

\subsection{Sampling}
	\label{sec:ppdpp-sampling}
The mixture representation provides the main steps of the sampling algorithm for extended L-ensembles. First of all, note that step 3 of both representations (the varying-size and the fixed-size cases) is a projection DPP, for which a generic algorithm is given in Algorithm~\ref{alg:proj_DPP} (see for instance~\cite{tremblay_optimized_2018}).
 \begin{algorithm}
 	\caption{Sampling $\mathcal{X}\sim |DPP|_m(\bU\bU^\top)$ with $\bU\in\mathbb{R}^{n\times m}$ verifying $\bU^\top \bU= \bI_m$.}
 	\label{alg:proj_DPP}
 	\begin{algorithmic}
 		\normalsize
 		\Input $\bU\in\mathbb{R}^{n\times m}$ verifying $\bU^\top \bU=\bI_m$\\
 		Write $\forall i, ~~\bm{y}_i=\bU^\top\bm{\delta_i}\in\mathbb{R}^{m}$.\\
 		$\mathcal{X} \leftarrow \emptyset$\\
 		Define $\bm{p}\in\mathbb{R}^n$ : $\forall i, \quad p(i) = \bm{y}_i^\top \bm{y}_i$\\
 		\textbf{for} $j=1,\ldots,m$ \textbf{do}:\\
 		\hspace{0.5cm}$\bm{\cdot}$ Draw $s_j$ with probability  $\mathbb{P}(s)=p(s)/\sum_{i}p(i)$\\
 		\hspace{0.5cm}$\bm{\cdot}$ $\mathcal{X} \leftarrow \mathcal{X}\cup\{s_j\}$\\
 		\hspace{0.5cm}$\bm{\cdot}$ Compute 
 		$\bm{f}_j = \bm{y}_{s_j} - \sum_{l=1}^{j-1} \bm{f}_l(\bm{f}_l^\top\bm{y}_{s_j})\in\mathbb{R}^{m}$\\
 		\hspace{0.5cm}$\bm{\cdot}$ Normalize $\bm{f}_j \leftarrow \bm{f}_j / \sqrt{\bm{f}_j^\top\bm{y}_{s_j}}$\\
 		\hspace{0.5cm}$\bm{\cdot}$ Update $\bm{p}$ : $\forall i\quad p(i) \leftarrow p(i) - (\bm{f}_j^\top\bm{y}_i)^2$\\
 		\textbf{end for}
 		\Output $\mathcal{X}$ of size $m$.
 	\end{algorithmic}
 \end{algorithm} 
\subsubsection{Sampling an extended L-ensemble based on a generic NNP} The following steps sample an extended L-ensemble based on a generic NNP $(\bL, \bV)$:
\begin{enumerate}
	\item Do a QR decomposition of $\bV$ to compute $\bQ$ 
	\item Compute and diagonalize $\widetilde{\bL}=(\bI-\bQ\bQ^\top)\bL(\bI-\bQ\bQ^\top)$ to obtain its truncated spectral decomposition 
	$\widetilde{\bU}\widetilde{\bm{\Lambda}}\widetilde{\bU}^\top$.
	\item Sample indices $\mathcal{Y}\sim DPP(\widetilde{\bm{\Lambda}})$
	\item Run Algorithm~\ref{alg:proj_DPP} with input $\bU=\left[\bQ, \widetilde{\bU}_{:,\mathcal{Y}}\right]\in\mathbb{R}^{n\times (p+|\mathcal{Y}|)}$
\end{enumerate}
Sampling a fixed-size extended L-ensemble is mainly the same: the sole difference is that step 3 becomes:
\begin{enumerate}
	\item[3'.] Sample indices $\mathcal{Y}\sim |DPP|_{m-p}(\widetilde{\bm{\Lambda}})$
\end{enumerate}
In terms of computation cost, step 1 requires $\mathcal{O}(np^2)$ number of operations, step 2 $\mathcal{O}(n^3)$, step 3 $\mathcal{O}(n)$, step 4 $\mathcal{O}(nm^2)$ in the fixed-size case and $\mathcal{O}(n(p+\mu)^2)$ in average in the varying-size case, where $\mu$ is the average number of samples of the diagonal DPP of step 3. The total sampling cost is thus dominated by the eigendecomposition step and boils down to $\mathcal{O}(n^3)$.
\subsubsection{Sampling an extended L-ensemble based on a low rank NNP} A first obvious acceleration of the previous scenario is when the columns of $\bV$ are already orthonormal. $\bQ$ is then set to $\bV$ and one starts the sampling algorithm directly at step 2. A more interesting scenario that arises in many practical applications is when $\bL$ is given in a low rank form $\bL=\bm{\Psi}\bm{\Psi}^\top$ with $\bm{\Psi}\in\mathbb{R}^{n\times r}$ and $r$ typically much smaller than $n$. In this case, one may circumvent the costly eigendecomposition in dimension $n$ of step 2 by a much more efficient singular value decomposition. In fact, step 2 becomes:
\begin{enumerate}
	\item[2'.] Compute and perform the SVD of  $\bm{\Psi}^\top(\bI-\bQ\bQ^\top)\in\mathbb{R}^{r\times n}$ to obtain its truncated SVD $\bm{\Psi}^\top(\bI-\bQ\bQ^\top)=\widetilde{\bW}\widetilde{\bm{\Sigma}}\widetilde{\bU}^\top$. Let $\widetilde{\bm{\Lambda}}=\widetilde{\bm{\Sigma}}^2$.
\end{enumerate}
Computing $\bm{\Psi}^\top(\bI-\bQ\bQ^\top)$ requires $\mathcal{O}(nrp)$ number of operations and computing its SVD $\mathcal{O}(nr^2)$. The total sampling time thus reduces to $\mathcal{O}(n(r^2+(p+\mu)^2))$ for varying-size extended L-ensembles and to $\mathcal{O}(n(r^2+m^2))$ for the fixed-size case.

\subsubsection{Approximate sampling using a Gibbs sampler}
\label{sec:gibbs}

An alternative to exact sampling is to use a Gibbs sampler. The advantage is
that no eigendecomposition is required, which is advantageous for extended L-ensembles where
$\bL$ is full rank and $n$ is large. Noting $\mu=\E(|X|)$, the cost of each
iteration in the sampler is $\O(\mu^3)$ in a naive implementation, and $\O(n)$
steps are required for mixing of the chain (see \cite{anari2020log} for recent results). The total cost
is then $\O(n\mu^3)$, but a careful implementation can bring this down to
$\O(n\mu^2)$, making it competitive with the exact sampler in the low-rank case.

The basic Gibbs sampler for (variable-size) DPPs can be derived by treating the
point process as a binary string $\vect{s}$ where $\vect{s}_i=1$ if $i \in \X$
and $\vect{s}_i=0$ otherwise. At each step, a random coordinate $i \in \{ 1
\ldots n\} $ is picked, and $s_i$ is flipped with probability
$p_{acc} = \frac{p(\vect{s}')}{p(\vect{s}')+p(\vect{s})}$, where $\vect{s}'$ is $\vect{s}$
with $s_i$ flipped. If $s_i=0$, the flip consists in adding an item to $\X$,
otherwise it consists in removing one, sometimes called ``up-down'' moves in the
literature. In the fixed-size case up-down moves are replaced with swaps, where
an item is taken out of $\X$ and another one is added. Note that  $p_{acc}$ only
depends on the ratio $\frac{p(\vect{s}')}{p(\vect{s})}$, and so can be computed
without knowledge of the normalisation constant of the extended L-ensemble.
The key to efficient implementation lies in taking advantage of rank-one or block updates
induced by the sampler. For instance, if item $i$ is added to the set $\X$, we
need to compute the ratio of
\[   \det
  \begin{pmatrix}
    \bL_{\X} & \bV_{\X,:} \\
    (\bV_{\X,:})^\top & \matr{0} 
  \end{pmatrix}
\]
  to
  \[   \det
    \begin{pmatrix}
      \bL_{\X'} & \bV_{\X',:} \\
      (\bV_{\X',:})^\top & \matr{0}
    \end{pmatrix}
  \]
  where $\X' = \X \cup i$. The latter matrix can be permuted to the bordered-matrix
  form: 
  \[
    \bA = 
    \begin{pmatrix}
      \bL_{\X} & \bV_{\X,:} & \bL_{i,\X} \\
      (\bV_{\X,:})^\top & \matr{0} & \bV_{i,;} \\
      \bL_{i,\X} & \bV_{i,:}^\top & L_{i,i}
    \end{pmatrix}
  \]
  and we may obtain the required ratio of determinants from Cramer's rule via
  $(\bA^{-1})_{i,i}$. Deleting an item works similarly, so that an efficient
  implementation of the Gibbs sampler requires only using block matrix inverses
  and some book-keeping. 

\section{Application: perturbative limits of L-ensembles}
\label{sec:ppDPP-as-limits}

As stated in the introduction, partial-projection DPPs arise as limits of certain L-ensembles, and we exhibit here one such limit: the L-ensemble based on the linear perturbation of a (low-rank)
positive semi-definite matrix; \emph{i.e.}, we consider L-ensembles based on matrices of the form:
\begin{equation}
  \label{eq:unbalanced-L-ensemble}
  \bL_\varepsilon \triangleq \varepsilon \bA + \bV \bV^\top
\end{equation}
where $\bA$ has full rank\footnote{The case where $\bA$ is not full rank can also be studied, but it is
	more burdensome and not much more informative} $n$, $\bV$ has full column rank $p < n$, and $\varepsilon$ is the parameter that will tend to $0$. 

Thus $\bL_\varepsilon$ defined in  \eqref{eq:unbalanced-L-ensemble} is a regular matrix pencil.  One should think about this
scenario as constructing a kernel as a sum of (a) a few important features
contained in $\bV\bV^{\top}$ and (b) a generic kernel in $\bA$.

\subsection{Limit of fixed-size L-ensembles based on \texorpdfstring{$\varepsilon \bA + \bV \bV^\top$}{}}
We begin with the more straightforward fixed-size case.
We seek the
limiting process $\Xs$ of $\X_\varepsilon \sim \mDPP{m}(\bL_\varepsilon)$ as $\flatlim$.
The following theorem establishes the limiting distribution  using asymptotic expansions of the  determinants. 

\begin{theorem}
  \label{thm:partial-proj-limit}
  Let $\Xe \sim \mDPP{m}(\bL_\varepsilon)$, with $\bL_\varepsilon$ as in Eq.~\eqref{eq:unbalanced-L-ensemble}. The limiting process is: 
\[
\Xe \rightarrow \Xs \sim 
\begin{cases}
\mDPP{m}(\bV\bV^{\top}),  & \text{if}~~m \le p \\
\mppDPP{m} \ELE{\bA}{\bV}, & \text{if}~~m > p. \\
\end{cases}  
\]
 \end{theorem}
\begin{proof}

First, we consider the case $m \le p$.
Note that the unnormalized probability mass function for the $L$-ensemble based on $\bL_\varepsilon$ is 
\[
f_{\varepsilon} (X)  = \det((\varepsilon \bA +  \bV\bV^\top)_{X}) =  \det( \varepsilon  \bA _{X} +  \bV_{X,:}( \bV_{X,:})^\top)=  \det( \bV_{X,:}(\bV_{X,:})^\top)  + \O(\varepsilon).
\]
Since $\rank \bV = p \ge m$, there exists a subset of rows $X_0$ such that 
\begin{equation}\label{eq:VVT_nonzero_minor}
\det( \bV_{X_0,:}(\bV_{X_0,:})^\top) \neq 0.
\end{equation}
so that the first term in the expansion is non-zero for at least some $X$. 
Therefore, we get that $\Xe  \rightarrow \mDPP{m}(\bV\bV^{\top})$.

The case $m > p$ is more delicate, as eq. \eqref{eq:VVT_nonzero_minor} no longer holds true, and we need to determine the  order of $\varepsilon$ in the expansion of $f_{\varepsilon} (X)$.
For this, we can invoke lemma~\ref{lem:det-coef-polynomial}  and remark~\ref{rem:det-coef-degree} to get
\begin{align*}
f_{\varepsilon} (X)  &   =\det( \varepsilon  \bA _{X} +  \bV_{X,:}( \bV_{X,:})^\top) = \varepsilon^{m} \det(  \bA _{X} + \varepsilon^{-1}   \bV_{X,:}( \bV_{X,:})^\top)\\
 &= 
 \varepsilon^{m} \left(\varepsilon^{-p}  (-1)^p\det
    \begin{pmatrix}
      \bA_{X} & \bV_{X,:} \\
     (\bV_{X,:})^\top & \matr{0} 
    \end{pmatrix} + \varepsilon^{-(p-1)} \ldots   \right)   \\ 
& = \varepsilon^{m-p} \left(  (-1)^p\det
    \begin{pmatrix}
      \bA_{X} & \bV_{X,:} \\
     (\bV_{X,:})^\top & \matr{0} 
    \end{pmatrix} +\O(\varepsilon) \right).
\end{align*}
In the limit, we get
 \begin{equation}
    \label{eq:extended-L-ensemble_bis}
    \Proba(\Xs = X) \propto  (-1)^p\det
    \begin{pmatrix}
      \bA_{X} & \bV_{X,:} \\
      (\bV_{X,:})^\top & \matr{0} 
    \end{pmatrix},
  \end{equation}
and hence $\Xe  \rightarrow \mppDPP{m} \ELE{\bA}{\bV}$.  
\end{proof}

\subsection{Limits of variable-size L-ensembles}
\label{sec:variable-size-pdpps-as-limits}

The variable-size version of the results requires a bit more care. In fixed-size
L-ensembles, the law of $\X$ is invariant to a rescaling of the positive semi-definite matrix it is based on: $\X \sim
\mDPP{m}(\bL)$ is equivalent to $\mDPP{m}(\alpha\bL)$ for any $\alpha > 0$. For
regular (variable-size) DPPs this is not true. That feature both enriches and
complicates a little the asymptotic analysis. 

\subsubsection{A trivial limit}
\label{sec:trivial-limit}

Let us start with a straightforward limit, namely $\Xe \sim
DPP(\bL_\varepsilon)$ based on the matrix pencil defined in  \eqref{eq:unbalanced-L-ensemble}.  

\begin{proposition}
  \label{prop:standard-dpp-limit}
  Let $\Xe \sim DPP(\varepsilon\bA + \bV\bV^\top)$. Then the limiting process $\Xs$ is $\Xs \sim
  DPP( \bV \bV^\top)$.
\end{proposition}
\begin{proof}
  Follows from pointwise convergence of the probability mass function, noting
  that $\det(\bI + \varepsilon\bA + \bV\bV^\top) = \det(\bI
  +\bV\bV^\top)(1+\O(\varepsilon))$.
  
  % \[ \bK_\varepsilon = \bL_\varepsilon ( \bL_\varepsilon + \bI)^{-1}\] with
  % $\bL_\varepsilon = \varepsilon\bA + \bV\bV^\top $
\end{proof}

The result is not very surprising. It has a noteworthy consequence, which is
that as $\flatlim$, the expected sample size will  be bounded by $p$ from above. %:  
%\[ \E(|\Xe|) = \sum_{i=1}^n \pi_i(\varepsilon) = \sum_{i=1}^p \frac{  \lambda_{i,0}}{1+ \lambda_{i,0}} +
%  \O(\varepsilon) \le p + \O(\varepsilon).\] 
If we wish to sample a larger number of points on average, then it appears that we
are out of luck.

\subsubsection{A more interesting limit}

We may instead look at a very similar limit: instead of taking $\bL_\varepsilon$, we will now take 
\[
\bL'_\varepsilon
= \varepsilon^{-1} \bL_\varepsilon = \bA + \varepsilon^{-1}\bV\bV^\top,
\]
which
carries the same intuition of giving more importance to $\bV\bV^\top$ than $\bA$.
\begin{proposition}
  \label{prop:standard-dpp-limit-partial}
  Let $\Xe \sim DPP(\bA +
  \varepsilon^{-1} \bV\bV^\top)$. Then the limiting process is $\Xs \sim
  DPP\ELE{\bA}{\bV}$.
\end{proposition}
\begin{proof}
  In appendix \ref{sec:limit-proof}. 
\end{proof}

 Importantly, the expected sample size of this rescaled L-ensemble allows for
 sample sizes larger than $p$.

 \begin{example}
   We return to the distribution of roots in a random forest, and show that it
   arises naturally as a limit. 
   The ``resistance distance'' in a graph (see for instance~\cite{klein1993resistance}) between two nodes $i$ and $j$ is defined as: 
   \[ D_{ij} =
     \mathcal{L}^\dag_{ii}+\mathcal{L}^\dag_{ii}-2\mathcal{L}^\dag_{ij} \]
   where $\gL^\dag$ is the pseudo-inverse of the Laplacian. The formula is
   analoguous to $\norm{\vect{x}_i - \vect{x}_j}^2 = \vect{x}^\top \vect{x} +
   \vect{y}^\top \vect{y} - 2 \vect{x}^\top \vect{y} $, and in that sense the
   resistance distance is actually a squared distance.  A graph kernel can be
   constructed from the resistance distance as $L_{ij}= \exp(-D_{ij})$,
   analoguously to the Gaussian kernel. If we parametrise the kernel as
   $L_{ij}= \frac{\exp(-\varepsilon D_{ij})}{\varepsilon}$, then by taking a
   Taylor series in $\varepsilon$ we
   have 
   \[ \bL_\varepsilon = \varepsilon^{-1} \ones \ones^\top -  \bD + \O(\varepsilon) \]
   Let $\X_\varepsilon \sim DPP(\bL_\varepsilon)$, a DPP on the nodes on the
   graph. Then using 
   prop. \ref{prop:standard-dpp-limit-partial} we find that the limit of $\X_\varepsilon$ is
   the roots process $\X_\star \sim DPP\ELE{\gL^\dag}{\ones}$, as discussed in Section~\ref{sec:roots-forests}.
 \end{example}
 
The interested reader on the topic of such flat limits of DPPs is referred to~\cite{barthelme_determinantal_2021}, in which we study flat limits of general DPPs (and not only flat limits of a matrix pencil), making extensive use of the extended L-ensemble formalism.
% goes to:
%\[ \E(|\Xs|) = \sum_i^n \pi_i=1 = p + \sum_{i=p+1}^{n} \frac{ \lambda_{i,1}}{1+ \lambda_{i,1}} \geq p, \]
%so the rescaled L-ensemble allows for a larger sample size.

\section{Application: constructing DPPs from Conditional Positive Definite kernels}
\label{sec:ppdpp-examples-cpdef}

The Gaussian kernel is the default choice in machine learning wherever a kernel
is needed, and it is tempting to use it to formulate a ``default'' DPP for
selecting points in $\R^d$. We pick a spatial scale $l$, and build an L-ensemble
with $L_{ij}= \gamma \exp(-\frac{1}{2l^2}\norm{\vect{x_i} - \vect{x_j}}^2)$, where the
expected size of $\X$ can be set using $\gamma$. This gives a repulsive point
process with good empirical properties, provided that the spatial scale $l$ is
set to a reasonable value \cite{tremblay2019determinantal}.

If one does not wish to have to find such a scale, there is still a way to
formulate a reasonable ``default DPP'', based on extended L-ensembles and CPD
kernels. We show how in this section, and give numerical results. 

\subsection{Partial projection DPPs and conditional positive definite functions}

An important generalisation of positive definite kernels is the notion of
conditional positive definite kernels (see for example
\cite{micchelli1986interpolation},\cite{wendland2004scattered}). Conditional positive
definite kernels generate conditionally positive definite matrices when
evaluated at a finite set of locations, just like positive definite kernels
generate positive definite matrices. We will show here that extended L-ensembles
let us construct DPPs based on conditional positive definite functions.

Conditional positive definite functions are the continuous counterpart of NNPs.
In the literature, they are often defined with respect to polynomial basis
functions. A \emph{monomial} in $\R^d$ is a function $\vect{x}^{\vect{\beta}} =
\prod_{i=1}^d x_i^{\beta_i}$. The degree of a monomial equals the sum of the
degrees, i.e. $\sum_i \beta_i = |\vect{\beta}|$. For instance, $\vect{x}^{(2,1)}=x_1^2x_2$ is a
monomial of degree $3$. A polynomial is a weighted sum of monomials, with degree equal to
the maximum degree of its components. The monomials of degree $\leq l$
are a basis for the polynomials of degree $\leq l$.

We may now state the most common definition of conditionally positive definite functions: 
\begin{definition}
A function $f:\R^d \longrightarrow \R$ is conditionally positive definite of order $\ell$ if and only if, for any $n\in \mathbb{N}$, any $X= (\vect{x}_1,\ldots,\vect{x}_n) \in (\R^d)^n$, any $\vect{\alpha}\in \R^n$ satisfying $\sum_i \alpha_i \vect{x}_i^{\vect{\beta}}=0$ for all multi-indices $\vect{\beta}$  s.t. $|\vect{\beta}|<\ell$, the quadratic form
\[
 \sum_{i,j} \alpha_i \alpha_j f(\vect{x}_i-\vect{x}_j)
 \]
 is non-negative.
\end{definition}

Suppose now that we introduce kernel matrices $\vect{L}_{X} =
[f(\vect{x} -\vect{y})]_{\vect{x} \in X,\vect{y} \in X }$, and the multivariate Vandermonde matrix
$\vect{V}_{< \ell}( X)$. A multivariate Vandermonde matrix consists in the
monomials evaluated at the locations in $X$, with the locations along the rows
and the monomials along the columns:
\[\vect{V}_{< \ell}( X) = [ \vect{x}^\beta ]_{\vect{x} \in X, |\beta|  \leq \ell -1}\]

For instance, in dimension $d=2$, the monomials of degree $\ell \leq 2$ have
exponents in:
\[ \{ \vect{\beta}\ |\ |\vect{\beta}| \leq 2\} = \left\{ (0,0),(0,1),(1,0),(1,1),(2,0),(0,2)  \right\}\]
 and the Vandermonde matrix $\vect{V}_{< \ell}( X)$ will consists
of the concatenation of $[\vect{x}^{(0,0)}]_{\vect{x} \in X}$, $[\vect{x}^{(0,1)}]_{\vect{x} \in X}$, etc.

Then, an equivalent definition is
\begin{definition}
A function $f:\R^d \longrightarrow \R$ is conditionally positive definite of order $\ell$ if and only if, for any $n\in \mathbb{N}$, any $X= (\vect{x}_1,\ldots,\vect{x}_n) \in (\R^d)^n$, the matrix $\vect{L}_{X}$ is conditionally positive definite with respect to $\vect{V}_{< \ell}(X)$.
\end{definition}
This extends the possible functions used to measure diversity in DPP sampling.
For instance, the function $-\norm{\vect{x}-\vect{y}}$ is conditionally positive definite of order $\ell=1$,
meaning that the extended L-ensemble formed from
\[ \bL =  \left[- \gamma\norm{\vect{x}-\vect{y}} \right]_{\vect{x} \in
    \Omega,\vect{y} \in \Omega} \]
and $\bV = \ones$ is a valid choice for any $\gamma > 0$. We can still think of $\bL$ as a similarity
matrix, since the closer are $\vect{x}$ and $\vect{y}$ to each other, the larger (due to the minus sign) the associated entry in $\bL$. As we will see in the
numerical results below, the corresponding DPP is indeed repulsive, and it has
no spatial-length scale parameter. The only thing that needs to be set is
$\gamma$, to control the expected size of $\X$. 

Numerical results show that the degree of repulsiveness can be increased by
using as a kernel a (non-even) positive power of the distance function. Letting $r =
\norm{\vect{x}-\vect{y}}$, we have that $\phi(r)=(-1)^{\lceil \beta/2 \rceil}  r^\beta ; \beta >0, \beta \not\in
2\mathbb{N}$  is a conditional positive function of order $\lceil \beta/2
\rceil$. This tells us that
\begin{equation}
  \label{eq:distance-ele}
  \bL =  \left[ \gamma (-1)^{\lceil \beta/2 \rceil} \norm{\vect{x}-\vect{y}}^\beta \right]_{\vect{x} \in
    \Omega,\vect{y} \in \Omega}
\end{equation}
\[  \]
and $\bV = \bV_{< \lceil \beta/2
  \rceil }(\Omega)$ is a valid extended L-ensemble for $\Omega \subset \R^d$. We recover the basic
case outlined above with $\beta = 1$, and we stress that $\beta$ cannot be
an even integer. 
In practice, repulsiveness increases with $\beta$, but again no spatial
length-scale parameter is required, just a choice of $\gamma$ for the expected
size. 

\subsection{Sampling}
\label{sec:sampling}

A DPP based on the conditionally positive definite kernel of equation \eqref{eq:distance-ele} can be sampled
exactly using algorithm \ref{alg:proj_DPP} after an eigendecomposition with cost
$\O(n^3)$. To bring down that cost, some approximation is required. One way is
to use a Gibbs sampler, as outlined in section \ref{sec:gibbs}, to bring 
the cost down to $\O(nk^2)$ where $k = E(|\X|)$. Another is to replace the exact
eigendecomposition with an approximation that retains only the dominant
eigenvectors, for instance via the Lanczos method or random projections. Note
that even in this case a QR decomposition of $\bV$ needs to be performed, and
the number of columns of $\bV$ increases with $\beta$, so that there is a
trade-off between repulsiveness and sampling cost. Another aspect of increasing
$\beta$ is that it increases the minimum size of $|\X|$, which is bounded below
by the number of columns in $\bV$. 

\subsection{Empirical results}
\label{sec:numerics}

For the numerical experiments we sampled a ground set $\Omega$ from a standard
Gaussian in $\R^2$, with size $n=800$. We formed extended L-ensembles as in eq.
\eqref{eq:distance-ele} with different values of $\beta$, setting $\gamma$ such
that $E(|\X|) = 28$. All computations are done exactly using a Julia toolkit
developed by the authors\footnote{available at \textcolor{blue}{https://gricad-gitlab.univ-grenoble-alpes.fr/barthesi/dpp.jl}}.

\begin{figure}
  \centering
  \includegraphics[width=10cm]{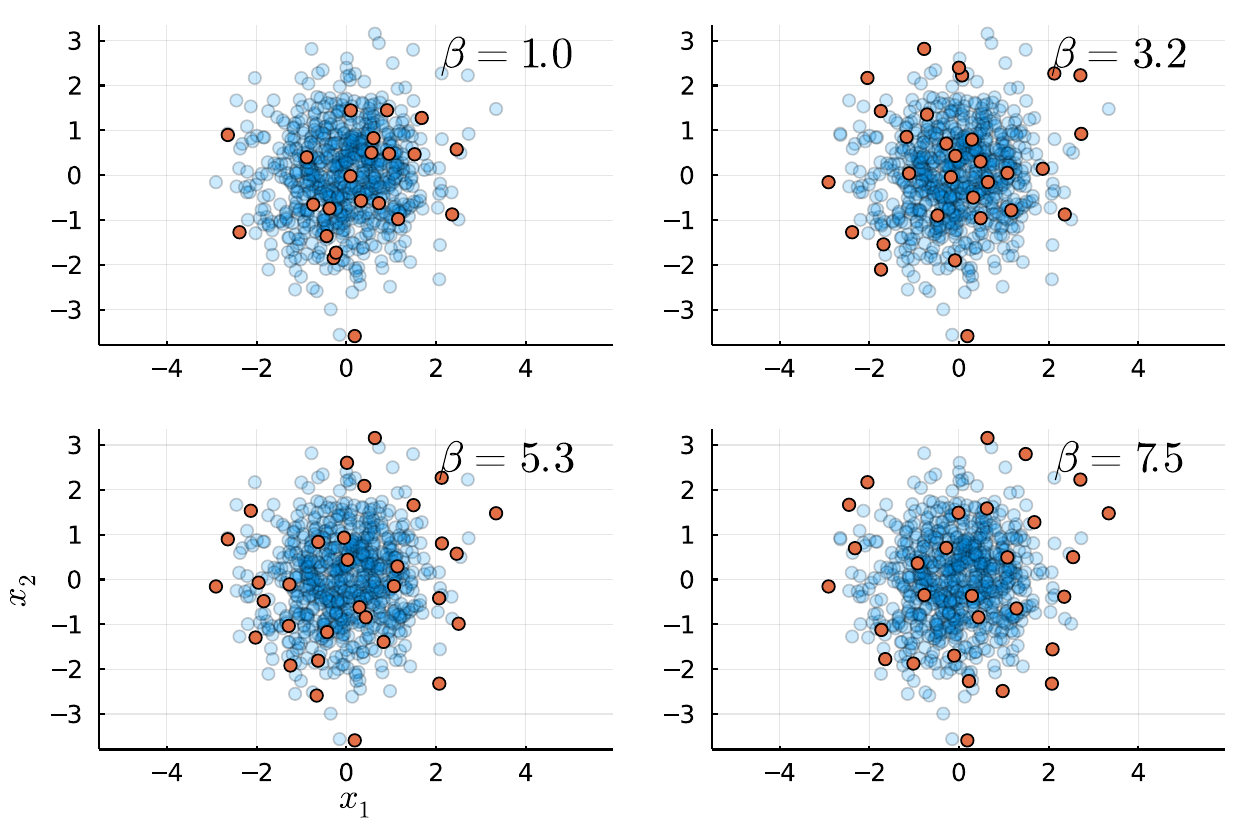}
  \caption{Samples from a DPP with ELE defined by eq. \eqref{eq:distance-ele}
    and different values of $\beta$.}
  \label{fig:samples}
\end{figure}

In figure \ref{fig:samples}, we show sampled points for different values of
$\beta$, with a visible increase in repulsiveness. This apparent increase can be
checked by examining the properties of the corresponding marginal kernel $\bK$.
One may use a repulsion index equivalent to the pair correlation function in
spatial statistics. We define:
\begin{equation}
  \label{eq:repulsion}
  \rho(i,j) = 1-\frac{p( i,j \subseteq \X )}{p(i \in \X)p(j \in \X)}
\end{equation}
which is 1 when the repulsion between $i$ and $j$ is very strong (they cannot be
sampled together) and 0 when the repulsion is null (they are included
``independently''). Using the definition of DPPs we obtain:
\[ \rho(i,j) = \frac{K_{ij}^2}{K_{ii}K_{jj}} \]

\begin{figure}
  \centering
  \includegraphics[width=10cm]{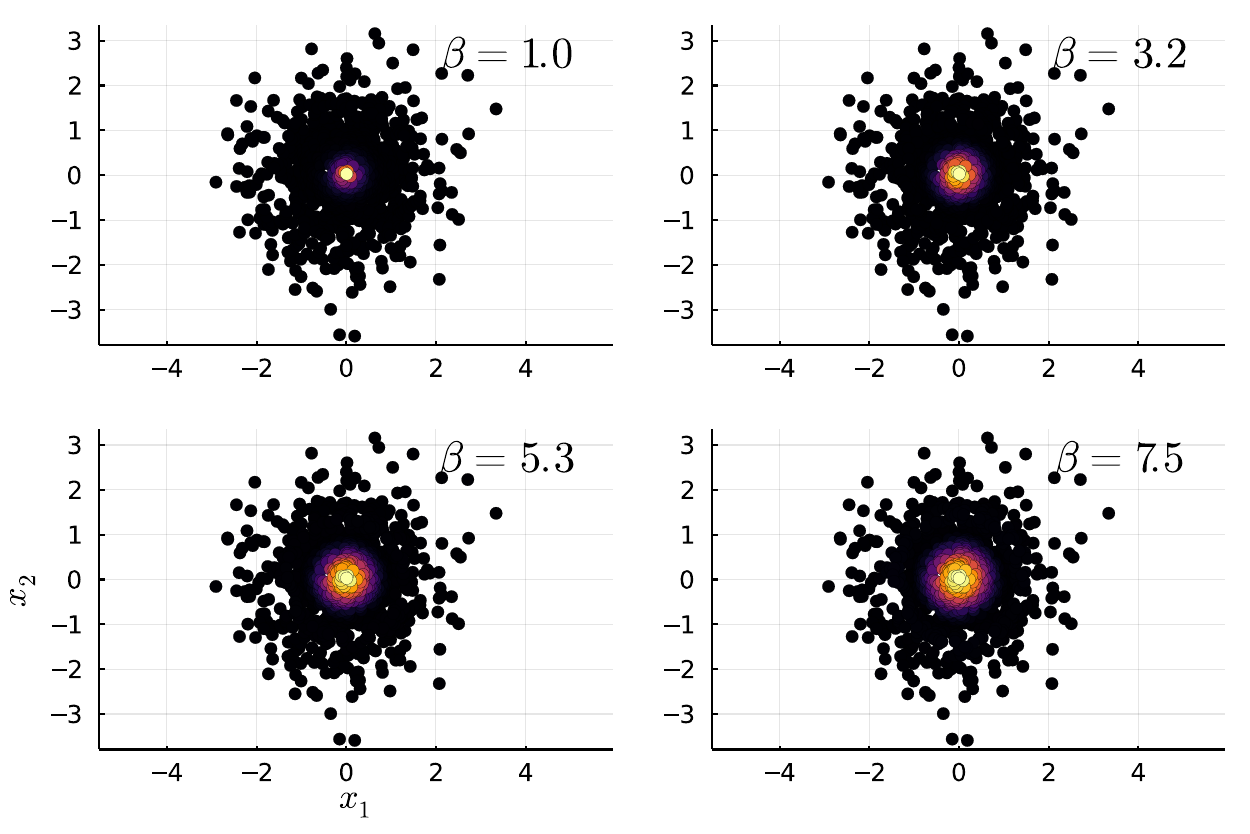}
  \caption{Repulsion index $\rho(i,j)$ for a fixed point $i$ near the center,
    and different values of $\beta$. The amount of repulsion is shown using a 
    colour scale, with white for the maximum value (1) and black for the
    minimum. The setup is the same as in fig. \ref{fig:samples}}
  \label{fig:repulsion}
\end{figure}

\begin{figure}
  \centering
  \includegraphics[width=10cm]{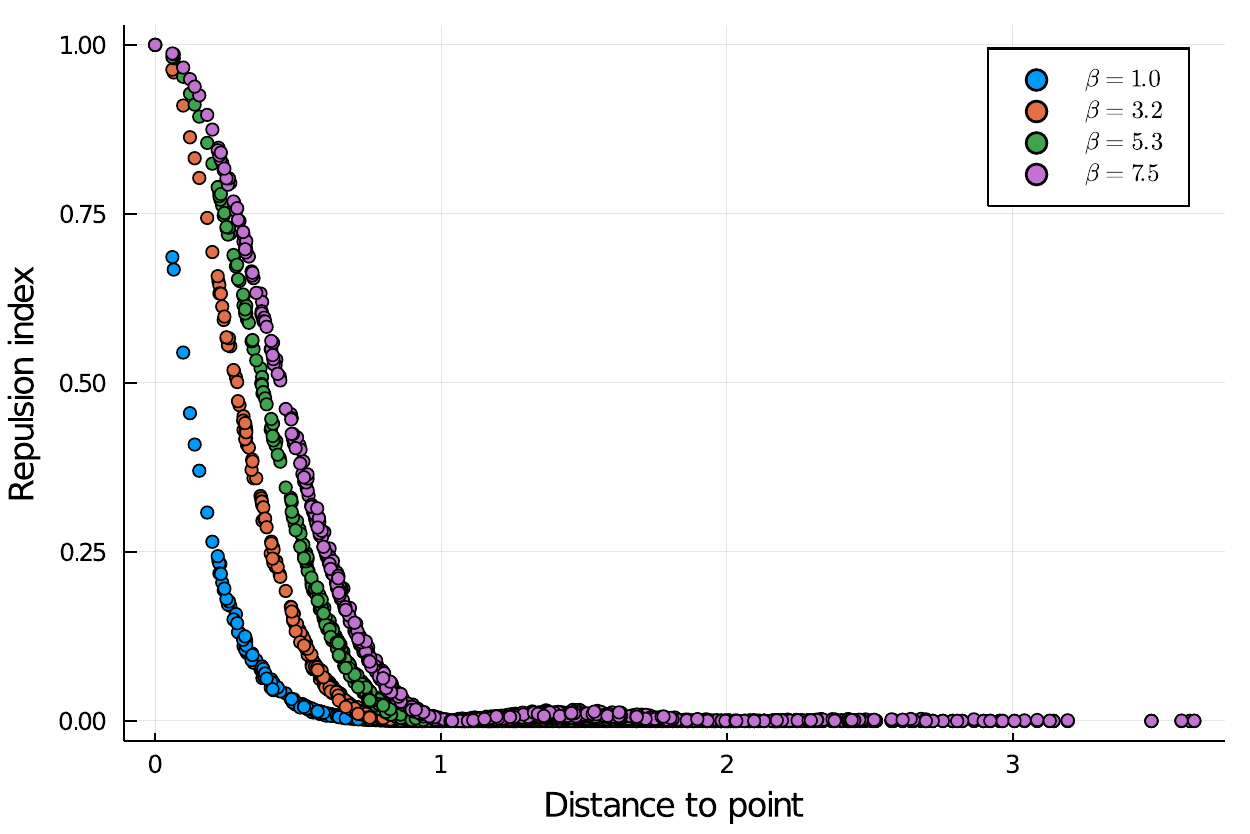}
  \caption{Repulsion index $\rho(i,j)$ for a fixed point $i$ near the center,
    plotted as a function of $\norm{\vect{x}_i - \vect{x}_j}$. The data are the same as in fig. \ref{fig:repulsion}}
  \label{fig:repulsion-dist}
\end{figure}

Figure \ref{fig:repulsion} shows $\rho$ plotted using a colour scale, with $i$ fixed and set to a
central point in $\Omega$. The increase in repulsiveness with $\beta$ is visible. Figure \ref{fig:repulsion-dist}
shows $\rho(i,j)$ plotted as a function of the distance $\norm{\vect{x}_i -
  \vect{x}_j}$. Since increasing $\beta$ increases the minimum size of $|\X|$,
we cannot increase $\beta$ indefinitely while keeping $\E |\X|$ constant. This
implies that there is a maximum amount of repulsiveness that can be achieved, in
line with known results on continuous DPPs \cite{biscio2016quantifying}.  

%As hinted in section  

\subsection{A link with interpolation}
\label{sec:link-interp}

There exists a link to interpolation (detailed recently
in~\cite{fanuel2020determinantal}) that we here illustrate. Suppose we want to
interpolate points $((\vect{x}_1,y_1), \ldots,(\vect{x}_n,y_n)\in (\R^d \times
\R)^n$ using the function
\[ s(\vect{x})=\sum_i \alpha_i f(\vect{x}-\vect{x}_i) +
  \sum_{k=1}^{\PP_{\ell-1,d}} \beta_k p_k(\vect{x}) \]
where $f$ is a conditionally positive function of order $\ell$, and $p_k$,
$k=1,\ldots, \PP_{\ell-1,d}$ is a basis for the set of polynomials of degree
less  or equal than $\ell-1$. For $s(x)$ to be an interpolant for the
measurements $y_1,y_2,\ldots,y_n$, we need to set $\vect{\alpha}$ and
$\vect{\eta}$ such that $s(\vect{x}_i)=y_i$ for all $i$. 
The solution of this interpolation problem is then equivalent to the solution of the linear system
\[
  \begin{pmatrix}
      \bL_{X} & \bV_{< \ell}\\
     (\bV_{< \ell})^\top & \matr{0} 
    \end{pmatrix}
      \begin{pmatrix}
      \vect{\alpha} \\
     \vect{\beta} 
    \end{pmatrix}
    =
      \begin{pmatrix}
     \vect{y}\\
    \vect{0} 
    \end{pmatrix}
  \]
where we recover the matrix defining the $L$-ensemble in partial projection
DPPs. A DPP based on the conditional positive definite kernel $f$ will sample a
good design for interpolation, since the interpolation points are
selected such that the interpolation matrix is well-conditioned. See
\cite{fanuel2020determinantal} for more on the topic.

\section{To conclude}
\label{sec:conclusion}

We have shown that DPPs can be represented as extended L-ensembles and that this
representation is both useful practically and theoretically. We have
deliberately kept to the discrete setting, in order to avoid the measure theory
that continuous point processes require. Extending our results to the continuous
setting is quite easy, and we would like to end this paper by sketching how.
Defining continuous extended L-ensembles is actually easier than the traditional construction of
DPPs from marginal kernels, since all that is needed is a conditionally positive definite function, to
define a valid likelihood (Janossi density). There is no need to verify that the
object exists, contrary to the standard approach. The next step is to use the
spectral theorem along with theorem \ref{thm:equivalence-extended-spectral}
(where the sum becomes infinite), to obtain that extended L-ensembles are mixtures of projection
DPPs, from which the properties of the marginal kernel follow. We may then use
the same ``default'' DPPs (eq. \eqref{eq:distance-ele}) introduced here in the
continuous setting. While the construction is quite straightforward, studying
the properties of these DPPs is much harder and an
interesting direction for future work.

\section*{Acknowledgments}
We thank Guillaume Gautier for helpful comments on preliminary versions of this
manuscript.

This work was supported by the ANR projects GenGP (ANR-16-CE23-0008), GRANOLA (ANR-21-CE48-0009), and 
LeaFleT (ANR-19-CE23-0021-01), as well as the LabEx PERSYVAL-Lab (ANR-11-LABX-0025-01),
the Grenoble Data Institute (ANR-15-IDEX- 02), MIAI@Grenoble Alpes
(ANR-19-P3IA-0003), the LIA CNRS/Melbourne Univ Geodesic, and the IRS (Initiatives de Recherche Stratégiques) of the IDEX Université Grenoble Alpes.

\begin{appendix}
  \label{sec:appendix}

\section{Inclusion probabilities in mixtures of projection DPPs}
\label{sec:marginal-kernel-ppDPPs-proof}

Here, we give formulas for inclusion probabilities valid for mixtures of projection
DPPs. These formulas yield the marginal kernels of L-ensembles and partial-projection
DPPs as a special case. We give a variant of a calculation in \cite{Barthelme:AsEqFixedSizeDPP},
appendix A.2.

Let $\bU$ be a fixed orthonormal basis of $\mathbb{R}^n$. We assume that $\X$ is generated according to the following mixture process:
\begin{enumerate}
\item Sample indices $\Y \sim \Proba(\Y)$
\item Form the projection matrix $\bM = \bU_{:,\Y} (\bU_{:,\Y})^\top$
\item Sample $\X | \Y \sim  \mDPP{m}(\bM)$
\end{enumerate}
We do not specify $\Proba(\Y)$ for now (it may be an L-ensemble, a fixed-size L-ensemble, etc.).

Since $\X$ is a mixture of projection-DPPs we can write 
\begin{eqnarray}
  \label{eq:incl-prob-mixture}
  \Proba(  \cW \subseteq \X)= \E_\Y[ \Proba( \cW \subseteq \X \vert \Y) ]
\end{eqnarray}
where the outer expectation is over $\Y$, the indices of the columns of $\bU$ 
sampled in the mixture process.
Since the innermost quantity is an inclusion probability for a projection DPP,
we have from lemma \ref{lem:marginal-kernel-proj}:
\begin{align*}
  \Proba( \cW \subseteq \X \vert \Y) & = \det \left(  \bM_\cW \right)\\
                                &=  \det \left( \bU_{\cW,\Y}(\bU_{\cW,\Y})^\top \right) \\
                                &=  \sum_{\cA \subseteq \Y,|\cA| = |\cW|}  \det\left( \bU_{\cW,\cA}\right)^2 
\end{align*}
where the last line follows from the Cauchy-Binet lemma (lemma
\ref{lem:cauchy-binet}). Injecting into \ref{eq:incl-prob-mixture}, we find:
\begin{align*}
  \Proba(  \cW \subseteq \X) &= \E_\Y\left[\sum_{\cA \subseteq \Y,|\cA| = |\cW|}  \det\left( \bU_{\cW,\cA}\right)^2 \right] \\
                        &= \sum_{\Y} \Proba(\Y) \sum_{\cA \subseteq \Y,|\cA| = |\cW|}  \det\left( \bU_{\cW,\cA}\right)^2 \\
                        &= \sum_{\Y,\cA \slash |\cA|=|\cW|} \det\left( \bU_{\cW,\cA}\right)^2 \Proba(\Y) \Ind\{\cA \subseteq \Y \} \\
                        &= \sum_{\cA \slash |\cA|=|\cW|} \det\left( \bU_{\cW,\cA}\right)^2 \Proba(\cA \subseteq \Y).
\end{align*}
In the case of L-ensembles and partial projection DPPs, we can go a bit further, since
the distribution of $\Y$ is a Bernoulli process (meaning that each element $i$
is included independently with probability $\pi_i$). In that case 
$\Proba(\cA \subseteq \Y) = \prod_{i \in \cA} \pi_i$, and using the Binet-Cauchy lemma 
once again we find: 
\begin{align*}
  \Proba(  \cW \subseteq \X) &= \sum_{\cA \slash |\cA|=|\cW|} \det\left( \bU_{\cW,\cA}\right)^2 \prod_{i \in \cA} \pi_i \\
                        &= \det \bU_{\cW,:} \diag(\pi_1,\ldots,\pi_n) (\bU_{\cW,:})^\top \\
                        &= \det \bK_{\cW} \numberthis \label{eq:marginal-kernel-generic}
\end{align*}
with $\bK = \bU \diag(\pi_1,\ldots,\pi_n) \bU^\top$.

  \section{Equivalence of extended L-ensembles and DPPs }
\label{sec:thm:K_to_L-ens_and_back}
We prove Th.~\ref{thm:L-ens_to_K} and \ref{thm:K_to_L-ens}.

\begin{proof}[{Proof of Th.~\ref{thm:L-ens_to_K} }]

	Let $\ELE{\bL}{\bV}$ be any NNP, and $\tbL$, $\bQ$, $\tbU$, $\tbLam$ and $q$ be as in Definition~\ref{def:nnp}. 
	Let $\X\in\Omega$ be drawn according to the distribution:
	\begin{align}
	\forall X\subseteq\Omega,\qquad \Proba(\X=X) \propto (-1)^p \det
	\begin{pmatrix}
	\bL_{X} & \bV_{X,:} \\
	(\bV_{X,:})^\top & \matr{0}
	\end{pmatrix}.
	\end{align} 
	Using the generalized Cauchy-Binet formula (theorem~\ref{thm:equivalence-extended-spectral}), this can be re-written as
	\begin{align}
	\forall X\subseteq\Omega,\qquad \Proba(\X=X) \propto \det(\bV^\top\bV)\sum_{Y,|Y|=m-p} \det \left( 
	\begin{bmatrix}
	\bQ_{X,:} & \tbU_{X,Y}
	\end{bmatrix}  \right)^2
	\prod_{i \in Y} \widetilde{\lambda}_i.
	\end{align} 
	As made precise by corollary~\ref{cor:mixture-rep-ppDPP-varying}, this equation can be interpreted from a mixture point of view. As such, the generic inclusion probability formulas of Appendix~\ref{sec:marginal-kernel-ppDPPs-proof} are applicable and one obtains the result.	
\end{proof}

\begin{proof}[{Proof of Th.~\ref{thm:K_to_L-ens}}]
	Given a marginal kernel $\matr{K}$, we can always rewrite  its spectral factorisation in the form of  Eq.~\eqref{eq:dpDPPmarginalkernel}, by grouping all the eigenvectors corresponding to the eigenvalue $1$ in $\matr{Q}$; all the remaining eigenvalues can be always represented as $\widetilde{\lambda}_i/(1+\widetilde{\lambda}_i)$.
\end{proof}

\section{Limit of perturbed L-ensembles}
\label{sec:limit-proof}

There are several ways to prove this result. The more elegant ones use standard
results from matrix perturbation theory. However, to avoid the introduction of additional
background, we prove the result by direct calculation, i.e. we
compute the limit of the marginal kernel. 

We seek to evaluate 
\[ \bK_\varepsilon = \bL_\varepsilon(\bL_\varepsilon + \bI)^{-1} = \bI - (\bL_\varepsilon + \bI)^{-1}\]
with $\bL_\varepsilon= \bA + \varepsilon^{-1} \bV\bV^\top$, which gives

\begin{equation}
  \label{eq:marginal-kernel-appendix}
  \bK_\varepsilon = \bI - \varepsilon (\varepsilon \bA + \varepsilon \bI+ \bV\bV^\top )^{-1}.
\end{equation}
\[  \]
We introduce a change of basis by the orthogonal matrix $\bR =
\begin{bmatrix}
  \bQ & \Qort
\end{bmatrix}
$
where as usual $\bQ$ spans the same column space as $\bV$. 
We apply the change of basis to $(\varepsilon \bA + \varepsilon \bI+ \bV\bV^\top
)$, to find:
\[ \varepsilon \bA + \varepsilon \bI+ \bV\bV^\top
   = \bR
   \begin{pmatrix}
     \varepsilon(\bP_1 + \bI) + \bT  & \varepsilon \bP_3 \\
     \varepsilon \bP_2 & \varepsilon (\bP_4 + \bI)
   \end{pmatrix} \bR^\top
\]
where $\bT = \bQ^\top \bV \bV^\top \bQ$, $\bP_1 = \bQ^\top \bA \bQ$, $\bP_2 = \Qort^\top \bA
\bQ$, $\bP_3 = \bP_2^\top$, $\bP_4 = \Qort \bA \Qort^\top$.
To find the inverse, we use the well-known formula for the inverse of block
matrices. The two Schur complements are as follows:
\[ \bS_a = \varepsilon(\bP_1 + \bI) + \bT + \varepsilon \bP_3 (\bP_4 + \bI)^{-1}
  \bP_2 = \bT + \O(\varepsilon) \]
and
\[ \bS_b = \varepsilon (\bP_4 + \bI) + \varepsilon^2 \bP_2 (\bT +
  \O(\varepsilon))^{-1} \bP_3 = \varepsilon ( \bP_4 + \bI  + \O(\varepsilon)) \]
Injecting into eq. \eqref{eq:marginal-kernel-appendix} and applying the
block-matrix formula, we find:
\[ \bR \bK_\varepsilon \bR^\top = \bI - 
  \begin{pmatrix}
    \matr{0} & \matr{0} \\
    \matr{0} & (\bP_4 + \bI)^{-1} 
  \end{pmatrix}
  + \O(\varepsilon)
\]
which implies:
\[ \bK_\varepsilon = \bQ\bQ^\top + \Qort (\Qort^\top \bA \Qort + \bI)^{-1} \Qort^\top =
  \bQ\bQ^\top + \tilde{\bA}(\tilde{\bA} + \bI)^{-1} \]
We may now apply th. \ref{thm:L-ens_to_K} and we are done. 
\end{appendix}

\bibliographystyle{imsart-number}
\bibliography{biblio.bib}

\end{document}